\documentclass{article}

\usepackage{amsmath}
\usepackage{amssymb}
\usepackage{cite}
\usepackage{tabularx}
\usepackage{graphicx}
\usepackage{caption}
\usepackage{subcaption}
\usepackage{authblk}
\usepackage{xcolor}
\usepackage{fontawesome5}
\usepackage{comment}

\title{\textbf{Convergence and almost sure exponential stability of compensated split-step theta scheme for stochastic pantograph models with Poisson random measure}}

\author{Amr Abosenna$^{1,*}$, Yongchun Zhou$^{2}$, Boping Tian$^{2}$\\
$^1${\textit{Department of Mathematics, Shanghai Normal University, Shanghai 200234, China}}\\
$^2${\textit{School of Mathematics, Harbin Institute of Technology, Harbin 150001, China}}\\
$^*${\textit{Corresponding authors: amr.abosena@feng.bu.edu.eg}}}\vspace{1em}
\date{}

\begin{document}
	\maketitle
	\graphicspath{{Pictures/}}
	\begin{abstract} 
		Recently, stochastic pantograph models have gained an intensive attention and have been used in different fields such as finance, biology, control and stochastic neural networks. It is also more preferable to incorporate jumps during the study of stochastic differential equations. In this paper, stochastic pantograph model with Poisson random measure is studied. The compensated split-step theta technique is applied to the considered model. The numerical scheme exhibits a non divergent attitude and converges to the solution of our model under assumptions addressed later on. Furthermore, the almost sure exponential stability of the numerical scheme is investigated via utilizing the discrete semi-martingale convergence theorem. Finally, theoretical findings are manifested via some numerical examples. 
	\end{abstract}

\section{Introduction}
Stochastic differential models are very important and many researchers have focused their attention on them because they have been widely used in many areas such as physics, chemistry, engineering, biology and mathematical finance to describe dynamical systems affected by uncertain factors. In order to gain more realistic simulations for stochastic systems, it is more desirable and efficient to study stochastic models with delay. In general, the future state of the stochastic system depends on the current and history states. Stochastic models accompanied with delay are named stochastic delay differential equations and perform better than stochastic models without delays \cite{babasola2023stochastic}. Hobson and Rogers \cite{hobson1998complete}, gave a new non-constant volatility model with past dependency in finance. Arriojas et al.\cite{arriojas2007delayed} assumed that the stock price follows a stochastic model with delay. Stochastic pantograph models \cite{vivek2022analysis,ren2023stability} are special kinds of stochastic delay differential equations with unlimited storage and are used in many fields of pure and applied mathematics such as probability and quantum mechanics. Ockendon and Tayler \cite{ockendon1971dynamics} studied the collection of the electric current via the pantograph of an electric locomotive, from which the name originates. In recent years, stochastic pantograph models have entered in many applications \cite{meng2011pathwise, appleby2016sufficient} and their importance appears in modelling random systems which have no instant impact from the time of their occurrence, but the feedback with delay has to be taken into consideration such as in control theory. Stochastic pantograph models with jumps are very important and can be applied in different fields of science and engineering. For example in financial Markets, Stochastic pantograph models are commonly used to describe the dynamics of financial markets, including stock prices, exchange rates, and interest rates allowing for better pricing and risk management in financial markets. In epidemic modelling, stochastic pantograph models can be employed to study the spread of infectious diseases and analyse the effectiveness of control strategies. These models can be used in simulating the epidemic's progression accurately, capturing the impact of delays and sudden changes in the infection rate, and aiding in designing effective intervention strategies. Furthermore in engineering systems, stochastic pantograph models find applications in various engineering systems, such as power systems, control systems, and communication networks. These models can also be used simulating the power grid's behaviour, accounting for sudden changes in demand, renewable energy generation, and other stochastic factors, thereby facilitating better system planning and control. \\

On the other hand, Weiner process is not the convenient approach for modelling situations having sudden changes and extreme events. Therefore, jump models are better for tackling these situations because they play a vital role in describing a sudden change of the system \cite{kou2002jump, svishchuk2000stochastic}. Merton \cite{merton1976option} was the first to propose jump-diffusion model. It is often better to use a jump-diffusion stochastic model for modelling stochastic systems \cite{tankov2003financial, bruti2007strong}. Stochastic models interspersed with Poisson jumps have been studied by many researchers \cite{maghsoodi1996mean,higham2005numerical,bruti2007approximation}. However, if the fluctuations are random process, then the number of points where jumps happen and the magnitude of these jumps are also stochastic and it is not enough to model these fluctuations using Poisson process and it is better to consider a general jump process arising from Poisson random measures and generated by Poisson point process. Stochastic models interspersed with L$\acute{e}$vy jumps have been investigated by many authors. Existence and uniqueness of stochastic models with Poisson random jumps have been examined by Bass et al. \cite{bass2004stochastic} and Albeverio et al. \cite{albeverio2010existence}. The asymptotic stability of stochastic models interspersed with L$\acute{e}$vy jumps was researched by  Applebaum \cite{applebaum2009asymptotic}. Furthermore, studying stochastic models with delay and jumps  is also preferable for better performance and accuracy \cite{agrawal2020jump,ahmadian2020exponential}. Therefore, this paper will deal with stochastic pantograph model with Poisson random measure.\\

However, most of stochastic pantograph models have difficulties in solving them analytically and numerical algorithms are needed to tackle this problem. Fan et al.\cite{fan2009alphath} presented numerical algorithms for solving stochastic pantograph differential equations via Razumikhin technique. Zhang et al. \cite{zhang2014convergence} focused their research on stochastic pantograph models and studied the convergence and stability of Euler-Maruyama technique. Fan et al. \cite{fan2007existence} studied the convergence of semi-implicit Euler techniques for stochastic pantograph models. Some authors focused the light on implicit techniques for stochastic differential models with and without Poisson jumps such as  \cite{hu2011convergence,buckwar2011runge ,hu1996semi}. Baker and Buckwar \cite{baker2000continuous} showed the existence of a unique solution of the linear stochastic pantograph model and they depicted that the numerical solution generated by the continuous theta algorithm converges to the exact solution with order half. However, all of the numerical schemes mentioned above are one-step explicit or implicit techniques.\\

In last years, a big concern was spotted on split-step theta methods because of their abilities in possessing the desirable stability characteristics and rates of convergence \cite{huang2014mean,liu2019split}. According to some authors \cite{liu2016mean, wang2011improved,zong2015exponential}, the split-step theta scheme has the ability of attaining the stability of the system with free choice of the stepsize. Furthermore, the inclusion of Euler-Maruyama scheme and split-step backward Euler scheme inside it by adjusting the value of theta to $ \theta=0 $ or $ \theta=1 $ respectively, can also be considered as an advantage. The split-step backward Euler technique applied to linear stochastic delay models was researched by Zhang et al. \cite{zhang2009split}. Cao et al. \cite{cao2014split} applied the split-step theta scheme to the stochastic delay models and showed that the numerical scheme has order of convergence of half and exponential mean square stability. Bao and Hu \cite{bao2023convergence} focused their research on the stochastic variable delay models and applied the split-step theta scheme on the model and studied its convergence and mean square stability.\\

However, as we know, there are not many literatures on the split-step theta schemes for stochastic pantograph models with Poisson random measure. Applying the stochastic theta scheme to stochastic pantograph model with Poisson random jumps and studying the convergence of it were done by Yang and Jiang \cite{yang2014stochastic}. Some publications put the light on presenting the convergence and stability of the numerical algorithms for stochastic models with Poisson random measure such as \cite{buckwar2011runge,hu2011convergence}. The convergence of Euler-Maruyama technique for a category of hybrid stochastic models incorporating L$\acute{e}$vy bursts was examined by Mao and Wei \cite{mao2013approximate}. Guo and Li \cite{guo2018almost} concentrated their research on studying the almost sure exponential stability of numerical solutions for stochastic pantograph models. Studying the stability of split-step theta scheme for neutral stochastic models with time lagging and Poisson leaps was done by Mo et al. \cite{mo2017exponential}. This paper will focus on examining the compensated split-step theta method's convergence and almost sure exponential stability for stochastic pantograph model augmented with Poisson random measure.\\

The following depicts how this paper is sorted. A collection of notations and the model description are shown in Section 2. Section 3 puts the light on the compensated split-step theta scheme and gives more elaboration about it. The convergence of the numerical algorithm is depicted in Section 4. Section 5 presents the almost sure exponential stability of the numerical scheme. Numerical examples are provided in Section 6 to foster the theoretical findings. Finally, some conclusions are mentioned in Section 7. 
	
\section{Model description}
In the paper,let $\left(\Omega,\mathcal{F},\mathcal{P}\right)$  be a complete probability space with filtration satisfying the usual conditions (\textit{i.e.}, it is increasing and right continuous while $\mathcal{F}_{0}$ contains all \textit{P}-null sets). Let $|\cdot|$ denotes the Euclidean vector norm in $  \mathbb R^{m} $ and let $\langle x,y\rangle$ be the inner product of $x$, $y$ in $\mathbb R^{m}$ and for $a \in\mathbb{R}$, $[a]$ refers to the integer part of $a$. Also, $a\vee b$ represents max($a,b$) and $a\wedge b$ represents min($a,b$). During our analysis, $ C $ is used to denote a general positive real constant and its value may differ at different positions.\\\\
Let $ W(t) $ be $ d $-dimensional Brownian motion defined on the same probability space and let $ U $ be a portion of $ \mathbb R^{m}\backslash \{0\} $ which is the range space of the burst jumps. Let $ N(\cdot,\cdot) $ defined on $ \mathbb R_{+}\times \mathbb R^{m} \backslash \{0\} $is  $\mathcal{F}_{t}$-adapted Poisson random measure which means that the process at any current time $t$ is known up to time $t$ and can not look forward into the future beyond time $t$. We have added that in the revised version and $ \widetilde{N}(dt,du)=N(dt,du)-\pi(du)dt, $ is the compensated Poisson random measure with L$\acute{e}$vy measure which is defined as the expected number of jumps of a certain size in a time interval of length one $\pi$ defined on $U$ with $\pi(U)=\lambda$. It is assumed that $ W(t) $ is independent of $ N(t,\cdot) $. For more details regarding the basics of stochastic analysis, see \cite{higham2001algorithmic}

Consider the following $m$-dimensional stochastic pantograph model with Poisson random measure of the form 
\begin{equation}
	\label{model}
	\begin{aligned}
		dx(t)&=f(x(t),x(\eta t))dt+g(x(t),x(\eta t))dW(t) \\
		&\quad+\int_{U}h(x(t),x(\eta t),u)N(dt,du),
	\end{aligned}
\end{equation}

defined on $ 0\leq t \leq T $ with  $0\leq\eta\leq 1 $ and initial data $x(0^-)=x_0$ where $x_0$ is $\mathcal{F}_{0}$-measurable, right-continuous and  $\mathbb{E}|x_{0}|^{p} < \infty $, for $p>0$. Let $ f(0,0)=g(0,0)=\int_{U}h(0,0,u)\pi(du)=0 $ which indicate that Eq.(\ref{model}) has a trivial solution. Here $f:\mathbb{R}^{m} \times \mathbb{R}^{m} \rightarrow \mathbb R^{m}$, $g:\mathbb R^{m}\times\mathbb R^{m}\rightarrow \mathbb R^{m\times d} $ and $h:\mathbb R^{m}\times\mathbb R^{m}\times U\rightarrow \mathbb R^{m} , m,d \in \mathbb{N}^{+} $ are Borel-measurable functions. In terms of compensated Poisson random measure, Eq.(\ref{model}) can be written as follows
\begin{equation}
	\label{compensated_model}
	\begin{aligned}
		dx(t)&=\bar{f}(x(t),x(\eta t))dt+g(x(t),x(\eta t))dW(t) \\
		&\quad+\int_{U}h(x(t),x(\eta t),u)\widetilde{N}(dt,du),
	\end{aligned}
\end{equation}
where $ \bar{f}(a,b)=f(a,b)+\int_U h(a,b,u)\pi(du)$ for all $ a, b \in \mathbb{R}^{m} $.\\\\
It is also convenient to write Eq.(\ref{compensated_model}) in stochastic integral form as follows
\begin{equation}
	\label{Integral_compensated_model}
	\begin{aligned}
		x(t)&=x_{0}+\int_{0}^{t}\bar{f}(x(s),x(\eta s))ds+\int_{0}^{t}g(x(s),x(\eta s))dW(s) \\
		&\quad+\int_{0}^{t}\int_{U}h(x(s),x(\eta s),u)\widetilde{N}(ds,du).
	\end{aligned}
\end{equation} 
\textbf{Assumption 2.1.} \textit{ For each $ b\geq1 $ and $ p\geq 2 $, there exists a positive constant $ k_{b} $ such that for all $ a_1,a_2,b_1,b_2 \in \mathbb{R}^{m} $ with $|a_1|\vee |a_2|\vee |b_1|\vee |b_2|\leq b $
\begin{equation}
	\label{Local_Lipschitz}
		|f(a_1,a_2)-f(b_1,b_2)|^{2}\vee|g(a_1,a_2)-g(b_1,b_2)|^{2}\leq k_{b}(|a_1-b_1|^{2}+|a_2-b_2|^{2})
\end{equation}
and
\begin{equation}
	\label{Local_Lipschitz_2}
		\int_{U}|h(a_1,a_2,u)-h(b_1,b_2,u)|^{p}\pi(du)\leq k_{b}(|a_1-b_1|^{p}+|a_2-b_2|^{p}).
\end{equation}}
 
\section{Compensated split-step theta method}
For simplification, we take equal partitions of the interval $[0,T]$ for a given large positive integer $\textit{K}$. Let $t_{n}=n\varDelta t$ for $n=0,1,\cdot\cdot\cdot,K$, where $\varDelta t =T/K $, then the compensated split-step theta scheme for Eq.(\ref{compensated_model}) is given by the initial value $Y_0=x_0$ and we compute $\{Y_n\}_{n=1}^{K}$ by
	\begin{equation}
	\label{scheme1}		
	Y_{n}^{*}=Y_n+ \theta \varDelta t \bar {f}(Y^{*}_{n},Y^{*}_{[\eta n]})
\end{equation}
\begin{equation}
	\label{scheme2}
	\begin{aligned}		
	Y_{n+1}&=Y_{n}+\varDelta t\bar{f}(Y^{*}_{n},Y^{*}_{[\eta n]})+g(Y^{*}_{n},Y^{*}_{[\eta n]})\varDelta W_n\\
	&\quad+\int_{t_n}^{t_{n+1}}\int_{U}h(Y^{*}_{n},Y^{*}_{[\eta n]},u)\widetilde{N}(dt,du)
\end{aligned}		
\end{equation}
for $n=0,1,\cdot\cdot\cdot,K-1$, where $\theta\in [0,1] $ is a fixed parameter, $Y_n$ approximates $x(t_n)$ at time $t_n$, $\varDelta W_n:=W(t_{n+1})-W(t_n)$ is the Brownian motion increment. The delay argument may not hit the previous time step which appears in the numerical method while dealing with the pantograph delay. This problem is tackled by interpolating the unknown approximate values of the solution to the closet grid point on the left endpoint of the interval containing the delay argument using piecewise constant polynomials \cite{zhang2014convergence}. For all $ t \in [t_{n},t_{n+1}) $, we define
\begin{equation}
	\label{continuous_approximation}
	\begin{aligned}	
		Y(t):&=Y_n+(t-t_n)\bar {f}(Y^{*}_{n},Y^{*}_{[\eta n]})+ g(Y^{*}_{n},Y^{*}_{[\eta n]})(W(t)-W(t_{n}))\\
		&\quad+\int_{t_n}^{t}\int_{U}h(Y^{*}_{n},Y^{*}_{[\eta n]},u)\widetilde{N}(ds,du)	
	\end{aligned}
\end{equation}
and denote 
\begin{equation}
	\begin{aligned}
	y(t)=\sum_{r=0}^{\infty}Y^{*}_{r}I_{ [t_{r},t_{r+1})}(t)
	\end{aligned}
\end{equation}
and
\begin{equation}
	\begin{aligned}
		z(t)=\sum_{r=0}^{\infty}Y^{*}_{[\eta r]}I_{ [t_{r},t_{r+1})}(t)
	\end{aligned}
\end{equation}
where $ I_{B} $ is the indicator function of set $ B $. Then, the continuous time approximation (\ref{continuous_approximation}) can be written in integral form as follows
\begin{equation}
	\label{SSTM_integral}
	\begin{aligned}
		Y(t)&=Y_{0}+\int_{0}^{t}\bar{f}(y(s),z(s))ds+\int_{0}^{t}g(y(s),z(s))dW(s)\\
		&\quad+\int_{0}^{t}\int_{U}h(y(s),z(s),u)\widetilde{N}(ds,du).
	\end{aligned}
\end{equation}
\textbf{Assumption 3.1.} \textit{The $p$th moments of the exact and continuous-time approximation solutions are bounded, that is, for some $ p\geq1$, there exists a positive constant $ C=C(p,T) $ such that 
\begin{equation}
	\mathbb{E}\left[ \sup_{0\leq t\leq T}\lvert x(t)\lvert^{p}\right]\bigvee \mathbb{E}\left[ \sup_{0\leq t\leq T}\lvert Y(t)\lvert^{p}\right]\leq C
\end{equation}}
To ensure that the compensated split-step theta scheme is well-defined, the following one-sided Lipschitz assumption \cite{wu2010almost} is imposed on $\bar{f}(p,q)$ in $p$.\\\\
\textbf{Assumption 3.2.} \textit{There exists a constant $ \kappa>0 $ such that for any  $ p_1 $, $ p_2 $, $ q \in \mathbb{R}^{m}$
	\begin{equation}
		\label{assump4.11}
		\begin{aligned}
			\langle p_1-p_2,\bar{f}(p_1,q)-\bar{f}(p_2,q)\rangle\leq \kappa|p_1-p_2|^2.
		\end{aligned}	
\end{equation}}
Under Assumption 3.2, if $\kappa\theta\Delta t<1$, then the numerical scheme is well-defined and has a unique solution \cite{higham2002strong}.

\section{Convergence of compensated split-step theta method}
To obtain our convergence result, we will first set up some helpful lemmas under the global Lipschitz condition and then make some relaxation and remove this condition by applying the stopping time technique.\\\\
\textbf{Assumption 4.1.} \textit{For any $ a_1, b_1 \in \mathbb{R}^{m} $ and there exists a positive constant $ K $ such that 
\begin{equation}
	\label{Linear_condition}
	\begin{aligned}
		|f(a_1,b_1)|^{2}\vee|g(a_1,b_1)|^{2}\leq K(1+|a_1|^{2}+|b_1|^{2})
	\end{aligned}
\end{equation}}
\textit{and
\begin{equation}
	\label{Linear_condition2}
	\begin{aligned}
		\int_{U}|h(a_1,b_1,u)|^{p}\pi(du)\leq K(1+|a_1|^{p}+|b_1|^{p}), \quad p\geq 2.
	\end{aligned}
\end{equation}}
\textbf{Lemma 4.1.} \textit{Suppose that Assumption 4.1 holds, if $ \varDelta t<\min \{1,1/C_p\} $, where $ C_p=C(p,\lambda,K,T,\theta)$ is a constant which depends on $ p, \lambda, K,T $ and $\theta $, then for every $ p\geq 2 $, there exists a positive constant $B_{1}=B_1(p,\lambda,K,T,\theta)$ independent of $\varDelta t$ such that\\
\begin{equation}
		\mathbb{E}\left[ \sup_{0\leq n\varDelta t\leq T}\lvert Y_{n}\lvert^{p}\right]\bigvee \mathbb{E}\left[ \sup_{0\leq n\varDelta t\leq T}\lvert Y^{*}_{n}\lvert^{p}\right]\leq B_{1} < \infty
\end{equation}
where $ Y_{n}^{*}$ and $ Y_{n} $ are defined in (\ref{scheme1}) and (\ref{scheme2}) respectively.}\\\\
\textbf{Proof.} By squaring both sides of (\ref{scheme1}), we get
\begin{equation}
	\label{S28}
	\begin{aligned}
		|Y_{n}^{*}|^{2}&=|Y_{n}|^{2}+\theta^2 \varDelta t^2|\bar{f}(Y^{*}_{n},Y^{*}_{[\eta n]})|^{2}+2\theta \varDelta t\langle Y_n ,\bar{f}(Y^{*}_{n},Y^{*}_{[\eta n]}) \rangle.
	\end{aligned}
\end{equation}
By utilizing inequalities $2\langle a,b\rangle \leq |a|^{2}+|b|^{2}$ and $ |a+b|^2\leq2(|a|^2+|b|^2) $, we can rewrite (\ref{S28}) as follows 
\begin{equation}
	\label{S3}
	\begin{aligned}
		|Y_{n}^{*}|^{2}&\leq(1+\theta \varDelta t)|Y_{n}|^{2}+2(\theta^2 \varDelta t^2+\theta \varDelta t)|f(Y^{*}_{n},Y^{*}_{[\eta n]})|^{2}\\
		&\quad+2(\theta^2 \varDelta t^2+\theta \varDelta t)\Biggl|\int_{U}h(Y^{*}_{n},Y^{*}_{[\eta n]},u)\pi(du)\Biggl|^{2}.\\
	\end{aligned}
\end{equation}
Raising both sides of (\ref{S3}) to the power $ p/2 $, then applying the inequality
\begin{equation}
	\label{bypass}
	\Biggl(\sum_{j=1}^{n}v_j\Biggl)^p\leq n^{p-1}\sum_{j=1}^{n}|v_j|^{p}
\end{equation}	
and taking the supremum and expectation, lead to
\begin{equation}
\label{help7}
\begin{aligned}
	\mathbb{E}\left[ \sup_{0\leq j\leq n}|Y_{j}^{*}|^{p}\right]
	&\leq3^{\frac{p}{2}-1}\Biggl\{(1+\theta\varDelta t)^{p/2}	\mathbb{E}\left[ \sup_{0\leq j\leq n}|Y_{j}|^{p}\right]\\
	&\quad+C_1\varDelta t^{p/2}\mathbb{E}\left[ \sup_{0\leq j\leq n}|f(Y^{*}_{j},Y^{*}_{[\eta j]})|^{p}\right]\\
	&\quad+C_1\varDelta t^{p/2}\mathbb{E}\left[ \sup_{0\leq j\leq n}\Biggl|\int_{U}h(Y^{*}_{j},Y^{*}_{[\eta j]},u)\pi(du)\Biggl|^{p}  \right] \Biggl\},
\end{aligned}		
\end{equation}
where $ C_1=(2(\theta^2+\theta))^{p/2}$. After applying Assumption 4.1 and noticing that $ \theta\in[0,1] $ and $ \varDelta t\leq1 $, we have 	
\begin{equation}
	\label{help8}
	\begin{aligned}
		\mathbb{E}\left[ \sup_{0\leq j\leq n}|Y_{j}^{*}|^{p}\right]
		&\leq3^{\frac{p}{2}-1}\Biggl\{(1+\theta\varDelta t)^{p/2}	\mathbb{E}\left[ \sup_{0\leq j\leq n}|Y_{j}|^{p}\right]\\
		&\quad+C_2\varDelta t^{p/2}\mathbb{E}\left[\sup_{0\leq j\leq n} (1+|Y_{j}^{*}|^{2}+|Y^{*}_{[\eta j]}|^{2})^{p/2} \right] \Biggl\},
	\end{aligned}		
\end{equation}
where $ C_2=K^{p/2}(1+\lambda^{p/2})C_1$.	
\begin{equation}
	\label{help9}
	\begin{aligned}
		\mathbb{E}\left[ \sup_{0\leq j\leq n}|Y_{j}^{*}|^{p}\right]
		&\leq C_3\mathbb{E}\left[ \sup_{0\leq j\leq n}|Y_{j}|^{p}\right]+C_4\varDelta t^{p/2}\\
		&\quad+2C_4\varDelta t^{p/2}\mathbb{E}\left[ \sup_{0\leq j\leq n}|Y^{*}_{j}|^{p} \right],
	\end{aligned}		
\end{equation}
where $ C_3=3^{\frac{p}{2}-1}2^{\frac{p}{2}}$ and $ C_4=3^{p-2}C_2$.
Therefore
\begin{equation}
	\label{S33}
	\begin{aligned}
		(1-2C_4\varDelta t^{p/2})\mathbb{E}\left[ \sup_{0\leq j\leq n}|Y^{*}_{j}|^{p} \right]&\leq C_3\mathbb{E}\left[ \sup_{0\leq j\leq n}|Y_{j}|^{p}\right]+C_4\varDelta t^{p/2}.
	\end{aligned}
\end{equation}
Because of $ \varDelta t<\min \{1,1/C_p\} $, where $ C_p=\frac{1}{(4C_4)^{2/p}}$, we have
\begin{equation}
	\label{help11}
	\begin{aligned}
		\mathbb{E}\left[ \sup_{0\leq j\leq n}|Y^{*}_{j}|^{p} \right]&\leq C\mathbb{E}\left[ \sup_{0\leq j\leq n}|Y_{j}|^{p}\right]+C
	\end{aligned}
\end{equation}
Let $ N $ and $ M $ be two positive integers satisfying $N\varDelta t\leq M\varDelta t\leq T$. By summing up (\ref{scheme2}) with respect to $\textit{n}$ from $j=0$ to $j=\textit{N-1}$, we get 
\begin{equation}
	\label{main1}
	\begin{aligned}
		 Y_{N} &= Y_{0}+ \varDelta t \sum_{j=0}^{\textit{N-1}}f(Y^{*}_{j},Y^{*}_{[\eta j]})
		 +\varDelta t\sum_{j=0}^{\textit{N-1}}\int_{U}h(Y^{*}_{j},Y^{*}_{[\eta j]},u)\pi(du)\\
		&\quad+\sum_{j=0}^{\textit{N-1}} g(Y^{*}_{j},Y^{*}_{[\eta j]})\Delta W_{j}
		+\sum_{j=0}^{\textit{N-1}}\int_{t_j}^{t_{j+1}}\int_{U}h(Y^{*}_{j},Y^{*}_{[\eta j]},u)\widetilde {N}(ds,du).
	\end{aligned}
\end{equation}
Raising both sides of (\ref{main1}) to the power $\textit{p}$, taking the supremum over $[0,\textit{M}]$ and mathematical expectation yield
\begin{equation}
	\label{term1}
	\begin{aligned}
	\mathbb{E}\left[\sup_{0\leq N\leq M}\lvert Y_{N}\lvert ^{p}\right] \leq 5^{p-1} \Biggl\{\mathbb{E}|Y_{0}|^{p}+T_1+T_2+T_3+T_4\Biggl\}, 
	\end{aligned}
\end{equation}
where
\begin{equation}
	T_1 := \varDelta t^{p}\mathbb{E}\left[ \sup_{0\leq N\leq M} \Biggl\lvert  \sum_{j=0}^{\textit{N-1}}f(Y^{*}_{j},Y^{*}_{[\eta j]})\Biggl\lvert^{p}\right] 
\end{equation}	

\begin{equation}
	T_2 := \varDelta t^{p}\mathbb{E}\left[ \sup_{0\leq N\leq M}\Biggl|\sum_{j=0}^{\textit{N-1}}\int_{U}h(Y^{*}_{j},Y^{*}_{[\eta j]},u)\pi(du)\Biggl|^{p}\right]
\end{equation}	

\begin{equation}
	T_3 := \mathbb{E}\left[ \sup_{0\leq N\leq M}\Biggl|\sum_{j=0}^{\textit{N-1}}g(Y^{*}_{j},Y^{*}_{[\eta j]})\Delta W_{j}\Biggl|^{p}\right] 
\end{equation}	
and
\begin{equation}
	T_4 := \mathbb{E}\left[ \sup_{0\leq N\leq M}\Biggl|\sum_{j=0}^{\textit{N-1}}\int_{t_j}^{t_{j+1}}\int_{U}h(Y^{*}_{j},Y^{*}_{[\eta j]},u)\widetilde {N}(ds,du)\Biggl|^{p}\right].
\end{equation}	
Next, the terms $ T_1, T_2, T_3$ and $ T_4 $ will be estimated separately. By utilizing (\ref{bypass}), Hölder inequality, Assumption 4.1 and getting use of (\ref{help11}), we have
\begin{equation}
	\label{term2}
	\begin{aligned}
		T_1
		&\leq \varDelta t^{p}M^{p-1}\mathbb{E}\sum_{j=0}^{\textit{M-1}}|f(Y^{*}_{j},Y^{*}_{[\eta j]})|^{p} \\
		&\leq \varDelta tT^{p-1}\mathbb{E}\sum_{j=0}^{\textit{M-1}}(|f(Y^{*}_{j},Y^{*}_{[\eta j]})|^{2})^{p/2}\\
		&\leq \varDelta tT^{p-1}K^{p/2}\mathbb{E}\sum_{j=0}^{\textit{M-1}}(1+|Y^{*}_{j}|^{2}+|Y^{*}_{[\eta j]}|^{2})^{p/2}\\
		&\leq C\varDelta t\sum_{j=0}^{\textit{M-1}}(1+\mathbb{E}|Y^{*}_{j}|^{p}+\mathbb{E}|Y^{*}_{[\eta j]}|^{p})\\
		&\leq C+C\varDelta t\sum_{j=0}^{\textit{M-1}}\mathbb{E}\left[\sup_{0\leq r\leq j} |Y^{*}_{r}|^{p}\right]\\
		&\leq C+C\varDelta t\sum_{j=0}^{\textit{M-1}}\mathbb{E}\left[\sup_{0\leq r\leq j} |Y_{r}|^{p}\right]
	\end{aligned}
\end{equation}
and
\begin{equation}
	\label{term4}
	\begin{aligned}
		T_2
		&\leq C\varDelta t^{p}M^{p-1}\mathbb{E}\sum_{j=0}^{\textit{M-1}}|\int_{U}h(Y^{*}_{j},Y^{*}_{[\eta j]},u)\pi(du)|^{p}\\
		&\leq C\varDelta tT^{p-1}\mathbb{E}\sum_{j=0}^{\textit{M-1}}\left[ \int_{U}\pi(du)\right] ^{p/2}\left[ \int_{U}|h(Y^{*}_{j},Y^{*}_{[\eta j]},u)|^{2}\pi(du)\right] ^{p/2}\\
		&\leq \varDelta tT^{p-1}[\pi(U)]^{p/2}K^{p/2}\mathbb{E}\sum_{j=0}^{\textit{M-1}}(1+|Y^{*}_{j}|^{2}+|Y^{*}_{[\eta j]}|^{2})^{p/2}\\
		&\leq C\varDelta t\sum_{j=0}^{\textit{M-1}}(1+\mathbb{E}|Y^{*}_{j}|^{p}+\mathbb{E}|Y^{*}_{[\eta j]}|^{p})\\
		&\leq C+C\varDelta t\sum_{j=0}^{\textit{M-1}}\mathbb{E}\left[\sup_{0\leq r\leq j} |Y^{*}_{r}|^{p}\right]\\
		&\leq C+C\varDelta t\sum_{j=0}^{\textit{M-1}}\mathbb{E}\left[\sup_{0\leq r\leq j} |Y_{r}|^{p}\right]. 
	\end{aligned}
\end{equation}
Then, by using discrete Burkholder-Davis-Gundy inequality \cite{mao2007stochastic} and proceeding as we did in (\ref{term2}), we have
\begin{equation}
	\label{term3}
	\begin{aligned}
		T_3
		&\leq D_p\mathbb{E}\left( \sum_{j=0}^{\textit{M-1}}|g(Y^{*}_{j},Y^{*}_{[\eta j]})|^{2}\varDelta t\right) ^{p/2}\\
		&\leq D_p\varDelta t^{p/2}M^{(p-2)/2}\mathbb{E}\sum_{j=0}^{\textit{M-1}}|g(Y^{*}_{j},Y^{*}_{[\eta j]})|^{p}\\
		&\leq D_p\varDelta tT^{(p-2)/2}K^{p/2}\mathbb{E}\sum_{j=0}^{\textit{M-1}}(1+|Y^{*}_{j}|^{2}+|Y^{*}_{[\eta j]}|^{2})^{p/2}\\
		&\leq C\varDelta t\sum_{j=0}^{\textit{M-1}}(1+\mathbb{E}|Y^{*}_{j}|^{p}+\mathbb{E}|Y^{*}_{[\eta j]}|^{p})\\
		&\leq C+C\varDelta t\sum_{j=0}^{\textit{M-1}}\mathbb{E}\left[\sup_{0\leq r\leq j} |Y^{*}_{r}|^{p}\right]\\
		&\leq C+C\varDelta t\sum_{j=0}^{\textit{M-1}}\mathbb{E}\left[\sup_{0\leq r\leq j} |Y_{r}|^{p}\right].
	\end{aligned}
\end{equation}
For the last term, by using the Kunita's first inequality \cite{kunita2004stochastic} and proceeding in a similar manner as we did before in (\ref{term2}) and (\ref{term4}), we get
\begin{equation}
	\label{term5}
	\begin{aligned}
		T_4
		&\leq \mu_p\varDelta t^{p/2}\mathbb{E}\left[ \sum_{j=0}^{\textit{M-1}}\int_{U}|h(Y^{*}_{j},Y^{*}_{[\eta j]},u)|^{2}\pi(du)\right]^{p/2}\\
		&\quad+\mu_p\varDelta t\mathbb{E}\sum_{j=0}^{\textit{M-1}}\int_{U}|h(Y^{*}_{j},Y^{*}_{[\eta j]},u)|^{p}\pi(du)\\
		&\leq C\varDelta t\sum_{j=0}^{\textit{M-1}}(1+\mathbb{E}|Y^{*}_{j}|^{p}+\mathbb{E}|Y^{*}_{[\eta j]}|^{p})\\
		&\leq C+C\varDelta t\sum_{j=0}^{\textit{M-1}}\mathbb{E}\left[\sup_{0\leq r\leq j} |Y^{*}_{r}|^{p}\right]\\
		&\leq C+C\varDelta t\sum_{j=0}^{\textit{M-1}}\mathbb{E}\left[\sup_{0\leq r\leq j} |Y_{r}|^{p}\right]. 
	\end{aligned}
\end{equation}
Substituting inequalities (\ref{term2}), (\ref{term4}), (\ref{term3}) and (\ref{term5}) into (\ref{term1}) yields 
\begin{equation}
	\label{final}
	\mathbb{E}\left[ \sup_{0\leq N\leq M}\lvert Y_{N}\lvert^{p}\right]\leq C+C \varDelta t\sum_{j=0}^{\textit{M-1}}\mathbb{E}\left[ \sup_{0\leq N\leq j}\lvert Y_{N}\lvert^{p}\right].
\end{equation}
Then, by applying the discrete-type Gronwall inequality to (\ref{final}), we obtain
\begin{equation}
	\mathbb{E}\left[ \sup_{0\leq n\varDelta t\leq T}\lvert Y_{n}\lvert^{p}\right]\leq C e^{CM\varDelta t}\leq B_{1}.
\end{equation}
Due to inequality (\ref{help11}), we get the boundedness of
\begin{equation}	
	\mathbb{E}\left[ \sup_{0\leq n\varDelta t\leq T}\lvert Y^{*}_{n}\lvert^{p}\right]\leq B_{1}.
\end{equation}
The proof is finished completely.\\\\
\textbf{Corollary 4.1.} \textit{Under Assumption 4.1 and for every $p\geq 2 $
\begin{equation}
	 \mathbb{E}\left[ \sup_{0\leq t\leq T}\lvert y(t)\lvert^{p}\right]\leq B_{2}<\infty,
\end{equation}}
where $ B_2=B_2(p,\lambda,K,T,\theta) $ is a positive constant which is independent of $ \varDelta t $.\\\\
\textbf{Proof.} Similar to the prove of Lemma 2.2 in \cite{jacob2009numerical}.\\\\
\textbf{Lemma 4.2.}	\textit{Under Assumption 4.1, if $ \varDelta t<\min \{1,1/C_p\} $, then for any $ p\geq2 $, there exists a positive constant $ C=C(p,\lambda,K,T,\theta) $ which is independent of $ \varDelta t $ such that 
\begin{equation}
	\label{lrr}
	\sup_{0\leq t\leq T}\mathbb{E}\Biggl[|Y(t)-y(t)|^{p}\Biggl]\leq C\varDelta t, \qquad \textit{for} \quad t \in [0,T]
\end{equation}}
\textbf{Proof}. Consider $t\in [t_n,t_{n+1})\subseteq [0,T]$. Then we have
\begin{equation}
	\label{lemma3.2}
	\begin{aligned}
		Y(t)&-y(t)\\
		&=Y(t)-Y_{n}\\
		&=(t-t_{n})f(Y^{*}_{n},Y^{*}_{[\eta n]})+(t-t_{n})\int_{U}h(Y^{*}_{n},Y^{*}_{[\eta n]},u)\pi(du)\\
		&\quad+g(Y^{*}_{n},Y^{*}_{[\eta n]})(W(t)-W(t_{n}))+\int_{t_{n}}^{t}\int_{U}h(Y^{*}_{n},Y^{*}_{[\eta n]},u) \widetilde{N}(ds,du).\\	
	\end{aligned}
\end{equation}
Raising both sides of (\ref{lemma3.2}) to the power $ p $ and taking expectation leads to
\begin{equation}
	\begin{aligned}
		\mathbb{E}|Y(t)-y(t)|^{p}
		&\leq 4^{p-1}\Biggl\{\mathbb{E}|(t-t_{n})f(Y^{*}_{n},Y^{*}_{[\eta n]})|^{p}\\
		&\quad+\mathbb{E}\Biggl|(t-t_{n})\int_{U}h(Y^{*}_{n},Y^{*}_{[\eta n]},u)\pi(du)\Biggl|^{p} \\
		&\quad+\mathbb{E}|g(Y^{*}_{n},Y^{*}_{[\eta n]})(W(t)-W(t_{n}))|^{p}\\
		&\quad+\mathbb{E}\Biggl|\int_{t_{n}}^{t}\int_{U}h(Y^{*}_{n},Y^{*}_{[\eta n]},u) \widetilde{N}(ds,du)\Biggl|^{p}\Biggl\}.
	\end{aligned}
\end{equation}
Upon using Hölder inequality, the property of martingale and seeking the same approach as in (\ref{term2}), (\ref{term4}), (\ref{term3}) and (\ref{term5}), we get
\begin{equation}
	\begin{aligned}
		\mathbb{E}|Y(t)-y(t)|^{p}
		&\leq 4^{p-1}\Biggl\{C\varDelta t^{p}\mathbb{E}\Biggl[|f(Y^{*}_{n},Y^{*}_{[\eta n]})|^{p}\Biggl]\\
		&\quad+C\varDelta t^{p}\mathbb{E}\Biggl[\int_{U}|h(Y^{*}_{n},Y^{*}_{[\eta n]},u)|^{p}\pi(du)\Biggl]\\
		&\quad+C\varDelta t^{\frac{p}{2}}\mathbb{E}\Biggl[|g(Y^{*}_{n},Y^{*}_{[\eta n]})|^{p}\Biggl]\\
		&\quad+C\varDelta t^{\frac{p}{2}}\mathbb{E}\Biggl[\int_{U}|h(Y^{*}_{n},Y^{*}_{[\eta n]},u)|^{2} \pi(du)\Biggl]^{\frac{p}{2}}\\
		&\quad+C\varDelta t\mathbb{E}\Biggl[\int_{U}|h(Y^{*}_{n},Y^{*}_{[\eta n]},u)|^{p} \pi(du)\Biggl]\Biggl\}.
	\end{aligned}
\end{equation}
Then by utilizing Assumption 4.1, we have 
\begin{equation}
	\begin{aligned}
	\mathbb{E}\Biggl[|Y(t)-y(t)|^{p}\Biggl]
	\leq 4^{p-1}C\varDelta t(1+\mathbb{E}|Y^{*}_{n}|^{p}+\mathbb{E}|Y^{*}_{[\eta n]}|^{p}).\\
	\end{aligned}
\end{equation}
By using Lemma 4.1, we get
\begin{equation}
	\label{lemma_3.20}
	\begin{aligned}
	\mathbb{E}\Biggl[|Y(t)-y(t)|^{p}\Biggl]\leq C\varDelta t,
	\end{aligned}
\end{equation}
where $ C $ is a positive constant which is independent of $ \varDelta t $. Therefore, we have
\begin{equation}
	\label{lemma_3.2dash}
	\begin{aligned}
	\sup_{0\leq t\leq T}\mathbb{E}\Biggl[|Y(t)-y(t)|^{p}\Biggl]\leq C\varDelta t.
	\end{aligned}
\end{equation}
The proof is complete.\\\\
\textbf{Lemma 4.3.}	\textit{Under Assumption 4.1, if $ \varDelta t<\min \{1,1/C_p\} $, then for any $ p\geq2 $, there exists a positive constant $ C=C(p,\lambda,K,T,\theta) $ which is independent of $ \varDelta t $ such that 
	\begin{equation}
		\label{lemma_002}
		\sup_{0\leq t\leq T}\mathbb{E}\Biggl[|Y(\eta t)-z(t)|^{p}\Biggl]\leq C\varDelta t, \qquad \textit{for} \quad t \in [0,T]	
	\end{equation}}
	\textbf{Proof}. Consider $t\in [t_n,t_{n+1})\subseteq [0,T]$. Then we have
	\begin{equation}
		\label{main_3.003}
		\begin{aligned}
			Y(\eta t)-z(t)&=Y(\eta t)-Y([\eta n]\varDelta t)\\
			&=\int_{[\eta n]\varDelta t}^{\eta t}\bar{f}(y(s),z(s))ds
			+\int_{[\eta n]\varDelta t}^{\eta t}g(y(s),z(s))dW(s)\\
			&\quad+\int_{[\eta n]\varDelta t}^{\eta t}\int_{U}h(y(s),z(s),u)\widetilde{N}(ds,du)\\
			&=\int_{[\eta n]\varDelta t}^{\eta t}f(y(s),z(s))ds+\int_{[\eta n]\varDelta t}^{\eta t}\int_{U}h(y(s),z(s),u)\pi(du)ds\\
			&\quad+\int_{[\eta n]\varDelta t}^{\eta t}g(y(s),z(s))dW(s)\\
			&\quad+\int_{[\eta n]\varDelta t}^{\eta t}\int_{U}h(y(s),z(s),u)\widetilde{N}(ds,du)\\	
		\end{aligned}
	\end{equation}
	Upon applying the inequality (\ref{bypass}) to (\ref{main_3.003}), we obtain
	\begin{equation}
		\label{origin}
		\begin{aligned}
			\mathbb{E}\Biggl[\sup_{[\eta n]\varDelta t\leq t\leq \eta t }|Y(\eta t)-z(t)|^{p}\Biggl]\leq 4^{p-1} \Biggl\{F_1+F_2+F_3+F_4\Biggl\}, 	
		\end{aligned}
	\end{equation}
	where
	\begin{equation}
		F_1:=\mathbb{E}\Biggl[\sup_{[\eta n]\varDelta t\leq t\leq \eta t }\Biggl|\int_{[\eta n]\varDelta t}^{t} f(y(s),z(s))ds\Biggl|^{p}\Biggl]	
	\end{equation}	
	\begin{equation}
		F_2:=\mathbb{E}\Biggl[\sup_{[\eta n]\varDelta t\leq t\leq \eta t }\Biggl|\int_{[\eta n]\varDelta t}^{ t}\int_{U} h(y(s),z(s),u)\pi(du)ds\Biggl|^{p}\Biggl]	
	\end{equation}
	\begin{equation}
		F_3:=\mathbb{E}\Biggl[\sup_{[\eta n]\varDelta t\leq t\leq \eta t }\Biggl|\int_{[\eta n]\varDelta t}^{ t} g(y(s),z(s))dW(s)\Biggl|^{p}\Biggl]
	\end{equation}
	and
	\begin{equation}
		F_4:=\mathbb{E}\Biggl[\sup_{[\eta n]\varDelta t\leq t\leq \eta t }\Biggl|\int_{[\eta n]\varDelta t}^{t}\int_{U} h(y(s),z(s),u)\widetilde{N}(ds,du)\Biggl|^{p}\Biggl].	
	\end{equation}
By applying Hölder inequality, Assumption 4.1, $ \eta t-[\eta n]\varDelta t\leq 2\varDelta t $ and Corollary 4.1, we get
\begin{equation}
	\label{first_term}
	\begin{aligned}
	F_1&\leq C\varDelta t^{p-1}	\mathbb{E}\int_{[\eta n]\varDelta t}^{\eta t} |f(y(s),z(s))|^{p}ds\\
	&\leq C\varDelta t^{p-1}\int_{[\eta n]\varDelta t}^{\eta t} (1+\mathbb{E}|y(s)|^{p}+\mathbb{E}|z(s)|^{p})ds\\
	&\leq C\varDelta t^p(1+2B_2)\\
	&\leq C\varDelta t^p	
	\end{aligned}
\end{equation}
and
\begin{equation}
	\label{third_term}
	\begin{aligned}
		F_2
		&\leq C\varDelta t^{p-1}\lambda\mathbb{E}\int_{[\eta n]\varDelta t}^{\eta t}\int_{U}|h(y(s),z(s),u)|^{p}\pi(du)ds\\
		&\leq C\varDelta t^{p-1}\int_{[\eta n]\varDelta t}^{\eta t} (1+\mathbb{E}|y(s)|^{p}+\mathbb{E}|z(s)|^{p})ds\\
		&\leq C\varDelta t^p(1+2B_2)\\
		&\leq C\varDelta t^p .
	\end{aligned}
\end{equation}
By using Burkholder-Davis-Gundy inequality \cite{mao2007stochastic}, Assumption 4.1 and Corollary 4.1, we have
\begin{equation}
	\label{second_term}
	\begin{aligned}
		F_3&\leq C\varDelta t^{\frac{p}{2}-1}\mathbb{E}\int_{[\eta n]\varDelta t}^{\eta t} |g(y(s),z(s))|^{p}ds\\
		&\leq C\varDelta t^{\frac{p}{2}-1}\int_{[\eta n]\varDelta t}^{\eta t} (1+\mathbb{E}|y(s)|^{p}+\mathbb{E}|z(s)|^{p})ds\\
		&\leq C\varDelta t^{\frac{p}{2}}(1+2B_2)\\
		&\leq C\varDelta t^{\frac{p}{2}}. 	
	\end{aligned}
\end{equation}
By using the Kunita's first inequality \cite{kunita2004stochastic}, Holder inequality, Assumption 4.1 and Corollary 4.1, we obtain
\begin{equation}
	\label{fourth_term}
	\begin{aligned}
		F_4
		&\leq C\mathbb{E}\Biggl[\int_{[\eta n]h}^{\eta t}\int_{U}|h(y(s),z(s),u)|^{2}\pi(du)ds\Biggl]^{\frac{p}{2}}\\
		&\quad+C\mathbb{E}\Biggl[\int_{[\eta n]\varDelta t}^{\eta t}\int_{U}|h(y(s),z(s),u)|^{p}\pi(du)ds\Biggl]\\
		&\leq (C\varDelta t^{\frac{p}{2}}+C)\int_{[\eta n]\varDelta t}^{\eta t} (1+\mathbb{E}|y(s)|^{p}+\mathbb{E}|z(s)|^{p})ds\\
		&\leq C\varDelta t(1+2B_2)\\
		&\leq C\varDelta t
	\end{aligned}
\end{equation}
Upon substituting (\ref{first_term}),(\ref{third_term}),(\ref{second_term}) and (\ref{fourth_term}) into (\ref{origin}), we conclude that
\begin{equation}
	\sup_{0\leq t\leq T}\mathbb{E}\Biggl[|Y(\eta t)-z(t)|^{p}\Biggl]\leq C\varDelta t,	
\end{equation}
where $ C $ is a positive constant which is independent of $ \varDelta t $. The proof is finished completely.\\\\
By using the stopping time technique, it is allowable to consider the solution of the model within a compact set where the coefficients are characterized to be Lipschitz. For each $ b>0 $, we define the stopping times\\
$ \nu_b=\inf\{t\in[0,T]:|x(t)|\geq b\}$, $ \vartheta_b=\inf\{t\in[0,T]:|Y(t)|\geq b\} $\quad and assume that $ \tau_b=\nu_b\wedge \vartheta_b$. Following the same approach as in Lemmas 4.2 and 4.3, the following corollaries can be obtained using the stopping time technique and proofs are omitted here.\\\\ 
\textbf{Corollary 4.2.}	\textit{Under Assumption 2.1, if $ \varDelta t<\min \{1,1/C_p^{b}\} $, where $ C_p^{b} $ is a constant depending on $b,p,\lambda,K_b,T$ and $\theta$. Then for any $ p\geq2 $, there exists a positive constant $ C_1(b) $ which is dependent on $ b $, but independent of $ \varDelta t $ such that 
	\begin{equation}
		\label{lemma_0002}
		\mathbb{E}|Y(t\wedge\tau_b)-y(t\wedge\tau_b)|^{p}\leq C_1(b)\varDelta t, \qquad \textit{for} \quad t \in [0,T].	
\end{equation}}
\textbf{Corollary 4.3.}	\textit{Under Assumption 2.1, if $ \varDelta t<\min \{1,1/C_p^{b}\} $, where $ C_p^{b} $ is a constant depending on $b,p,\lambda,K_b,T$ and $\theta$. Then for any $ p\geq2 $, there exists a positive constant  $ C_2(b) $ which is dependent on $ b $, but independent of $ \varDelta t $ such that 
	\begin{equation}
		\label{lemma_00002}
		\mathbb{E}|Y(\eta( t\wedge\tau_b))-z(t\wedge\tau_b)|^{p}\leq C_2(b)\varDelta t, \qquad \textit{for} \quad t \in [0,T].	
\end{equation}}
\textbf{Theorem 4.1.} \textit{Under Assumption 2.1, for any $ p\geq2,$ and $ \theta\in[0,1] $, the compensated split-step theta method approximate solution $ Y(t) $ converges in the $ p$th moment to the exact solution $ x(t) $, that is 
\begin{equation}
	\label{convergence}
	\lim_{ \varDelta t\rightarrow 0}\mathbb{E}\left[\sup_{0\leq t\leq T}|x(t)-Y(t)|^{p}\right] =0
\end{equation}}	
\textbf{Proof}. Let $ e(t)=x(t)-Y(t).$ Applying the Young inequality \cite{higham2002strong}, for any $ \rho>0 $ and $ \varrho>p,$ we have
\begin{equation}
\label{young1}
\begin{aligned}	
\mathbb{E}\Biggl[\sup_{0\leq t\leq T}|e(t)|^{p}\Biggl]&=\mathbb{E}\Biggl[\sup_{0\leq t\leq T}|e(t)|^{p}I_{\{\nu_b>T,\vartheta_b>T\}}\Biggl]\\
&\quad+\mathbb{E}\Biggl[\sup_{0\leq t\leq T}|e(t)|^{p}I_{\{\nu_b\leq T\,\text{or}\,\vartheta_b\leq T\}}\Biggl]\\
&\leq \mathbb{E}\Biggl[\sup_{0\leq t\leq T}|e(t)|^{p}I_{\{\tau_b>T\}}\Biggl]+\frac{\rho p}{\varrho}\mathbb{E}\Biggl[\sup_{0\leq t\leq T}|e(t)|^{\varrho}\Biggl]\\
&\quad+\frac{\varrho-p}{\varrho\rho^{\frac{p}{\varrho-p}}}\mathbb{P}(\nu_b\leq T\text{or}\,\vartheta_b\leq T).\\
\end{aligned}	
\end{equation}
Note that
\begin{equation}
\label{hhh1}	
\mathbb{E}\Biggl[\sup_{0\leq t\leq T}|e(t)|^{\varrho}\Biggl]\leq 2^{p-1}\mathbb{E}\Biggl[\sup_{0\leq t\leq T}(|x(t)|^{\varrho}+|Y(t)|^{\varrho})\Biggl]\leq 2^pC	
\end{equation}
and
\begin{equation}
\label{hhh2}
\begin{aligned}
\mathbb{P}(\vartheta_b\leq T)=\mathbb{E}\Biggl[I_{\{\vartheta_b>T\}}\frac{|Y(\vartheta_b)|^p}{b^p}\Biggl]\leq \frac{1}{b^p}\mathbb{E}\Biggl[\sup_{0\leq t\leq T} |Y(t)|^p\Biggl]\leq \frac{C}{b^p}.
\end{aligned}	
\end{equation}	
Same result can be obtained for $ \nu_b $. Thus
\begin{equation}
	\label{hhh3}
	\mathbb{P}(\nu_b\leq T\text{or}\,\vartheta_b\leq T)\leq\mathbb{P}(\vartheta_b\leq T)+\mathbb{P}(\nu_b\leq T)\leq\frac{2C}{b^p}.
\end{equation}	
By plugging in (\ref{hhh1}) and (\ref{hhh3}) into (\ref{young1}), we get 
\begin{equation}
	\label{young}
	\begin{aligned}	
		\mathbb{E}[\sup_{0\leq t\leq T}|e(t)|^{p}]&
		&\leq\mathbb{E}[\sup_{0\leq t\leq T}|e(t\wedge\tau_b)|^{p}]+2^{p}\frac{\rho p}{\varrho}C+\frac{\varrho-p}{\varrho\rho^{\frac{p}{\varrho-p}}}\frac{2C}{b^p}.
	\end{aligned}	
\end{equation}		
For the first term of (\ref{young}), we have
\begin{equation}
	\label{convergence1}
	\begin{aligned}
	\mathbb{E}\left[\sup_{0\leq t\leq T}|e(t\wedge\tau_b)|^{p}\right]
	\leq 4^{\frac{p-2}{2}}(G_1+G_2+G_3+G_4)
	\end{aligned}
	\end{equation}
where
\begin{equation*}
	G_1:=\mathbb{E}\Biggl[\sup_{0\leq t\leq T}\Biggl|\int_{0}^{t\wedge\tau_b}(f(x(s),x(\eta s))-f(y(s),z(s)))ds\Biggl|^{p}\Biggl]
\end{equation*}	
\begin{equation*}	
	G_2:=\mathbb{E}\Biggl[\sup_{0\leq t\leq T}\Biggl|\int_{0}^{t\wedge\tau_b}\int_{U}(h(x(s),x(\eta s),u)-h(y(s),z(s),u))\pi(du)ds\Biggl|^{p}\Biggl]
\end{equation*}	
\begin{equation*}
	G_3:=\mathbb{E}\Biggl[\sup_{0\leq t\leq T}\Biggl|\int_{0}^{t\wedge\tau_b}(g(x(s),x(\eta s))-g(y(s),z(s)))dW(s)\Biggl|^{p}\Biggl]
\end{equation*}
and
\begin{equation*}
	G_4:=\mathbb{E}\Biggl[\sup_{0\leq t\leq T}\Biggl|\int_{0}^{t\wedge\tau_b}\int_{U}(h(x(s),x(\eta s),u)-h(y(s),z(s),u))\widetilde{N}(ds,du)\Biggl|^{p}\Biggl].
\end{equation*}	
By applying Hölder inequality, Assumption 2.1, Corollaries 4.2 and 4.3, we get
\begin{equation}
	\label{conv_first_term}
	\begin{aligned}
	G_1
	&\leq C\mathbb{E}\left[\int_{0}^{T\wedge\tau_b}|f(x(s),x(\eta s))-f(y(s),z(s))|^{p}ds\right]\\
	&\leq C\mathbb{E}\left[\int_{0}^{T\wedge\tau_b}[k_l(|x(s)-y(s)|^{2}+|x(\eta s)-z(s)|^{2})]^{\frac{p}{2}}ds \right]\\
	&\leq C\mathbb{E}\left[\int_{0}^{T\wedge\tau_b}(|x(s)-Y(s)|^{p}+|Y(s)-y(s)|^{p})ds\right]\\
	&\quad+C\mathbb{E}\left[\int_{0}^{T\wedge\tau_b}(|x(\eta s)-Y(\eta s)|^{p}+|Y(\eta s)-z(s)|^{p})ds\right]\\
	&\leq C\mathbb{E}\left[\int_{0}^{T}(|x(s\wedge\tau_b)-Y(s\wedge\tau_b)|^{p}+|Y(s\wedge\tau_b)-y(s\wedge\tau_b)|^{p})ds\right]\\
	&\quad+C\mathbb{E}\left[\int_{0}^{T}(|x(\eta(s\wedge\tau_b))-Y(\eta(s\wedge\tau_b))|^{p}+|Y(\eta(s\wedge\tau_b))-z(s\wedge\tau_b)|^{p})ds\right]\\
	&\leq C\varDelta t+C\int_{0}^{T}\mathbb{E}\left[\sup_{0\leq r\leq s}|e(r\wedge\tau_b)|^{p}\right]ds
	\end{aligned}
\end{equation}
and
\begin{equation}
	\label{conv_second_term}
	\begin{aligned}
		G_2
		&\leq C\mathbb{E}\left[\int_{0}^{T\wedge\tau_b}\Biggl|\int_{U}(h(x(s),x(\eta s),u)-h(y(s),z(s),u))\pi(du)\Biggl|^{p}ds\right]\\
		&\leq C\mathbb{E}\int_{0}^{T\wedge\tau_b}\Biggl[\int_{U}\pi(du)\Biggl]^{\frac{p}{2}}\Biggl[|h(x(s),x(\eta s),u)-h(y(s),z(s),u)|^{2}\pi(du)\Biggl]^{\frac{p}{2}}ds\\
		&\leq C\mathbb{E}\int_{0}^{T\wedge\tau_b}\Biggl[|h(x(s),x(\eta s),u)-h(y(s),z(s),u)|^{2}\pi(du)\Biggl]^{\frac{p}{2}}ds.\\
	\end{aligned}
\end{equation}
Then by proceeding as we did before in (\ref{conv_first_term}), we obtain
\begin{equation}
\label{second2_con.}
\begin{aligned}
	G_2&\leq C\varDelta t+C\int_{0}^{T}\mathbb{E}\left[\sup_{0\leq r\leq s}|e(r\wedge\tau_b)|^{p}\right]ds.
\end{aligned}		
\end{equation}
Now, by utilizing the Burkholder-Davis-Gundy inequality \cite{mao2007stochastic} and Holder inequality, we have a positive constant $D_p $ such that	
\begin{equation}
	\label{conv_third_term}
	\begin{aligned}
		G_3
		&\leq D_p\mathbb{E}\left[\int_{0}^{T}|g(x(s\wedge\tau_b),x(\eta (s\wedge\tau_b)))-g(y(s\wedge\tau_b),z(s\wedge\tau_b))|^{2}ds\right]^{\frac{p}{2}}\\
		&\leq D_pT^{\frac{p}{2}-1}\mathbb{E}\int_{0}^{T}|g(x(s\wedge\tau_b),x(\eta (s\wedge\tau_b)))-g(y(s\wedge\tau_b),z(s\wedge\tau_b))|^{p}ds.\\
	\end{aligned}
\end{equation}
Then, following the same procedures as in (\ref{conv_first_term}), we get
\begin{equation}
	\label{third3_con.}
	\begin{aligned}
		G_2&\leq C\varDelta t+C\int_{0}^{T}\mathbb{E}\left[\sup_{0\leq r\leq s}|e(r\wedge\tau_b)|^{p}\right]ds.
	\end{aligned}		
\end{equation}
By using the Kunita's first inequality \cite{kunita2004stochastic}, we have a positive constant $\mu_p $ such that
\begin{equation}
	\label{conv_fourth_term}
	\begin{aligned}
		G_4
		&\leq \mu_p\Biggl\{\mathbb{E}\left[\int_{0}^{T\wedge\tau_b}\int_{U}|h(x(s),x(\eta s),u)-h(y(s),z(s),u)|^{2}\pi(du)ds\right]^{\frac{p}{2}} \\
		&\quad+\mathbb{E}\left[\int_{0}^{T\wedge\tau_b}\int_{U}|h(x(s),x(\eta s),u)-h(y(s),z(s),u)|^{p}\pi(du)ds\right]\Biggl\}. 
	\end{aligned}
\end{equation}
Then, in a same manner as $ G_1 $ and $ G_2 $ were done, we can obtain
\begin{equation}
	\label{forth4_con.}
	\begin{aligned}
		G_2&\leq C\varDelta t+C\int_{0}^{T}\mathbb{E}\left[\sup_{0\leq r\leq s}|e(r\wedge\tau_b)|^{p}\right]ds.
	\end{aligned}		
\end{equation}
By plugging in (\ref{conv_first_term}), (\ref{second2_con.}), (\ref{third3_con.}) and (\ref{forth4_con.}) into (\ref{convergence1}), we get
\begin{equation}
	\label{convergence_final}
	\begin{aligned}
		\mathbb{E}\left[\sup_{0\leq t\leq T}|e(t\wedge\tau_b)|^{p}\right]
		\leq C\varDelta t+C\int_{0}^{T}\mathbb{E}\left[\sup_{0\leq r\leq s}|e(r\wedge\tau_b)|^{p}\right]ds.
	\end{aligned}
\end{equation}
Due to Gronwall inequality, we have 
\begin{equation}
	\label{convergence_final_1}
	\begin{aligned}
		\mathbb{E}\left[\sup_{0\leq t\leq T}|e(t\wedge\tau_b)|^{p}\right]
		\leq C\varDelta te^{CT}.
	\end{aligned}
\end{equation}
By plugging in (\ref{convergence_final_1}) into (\ref{young}), we have
\begin{equation}
\label{final10}
\begin{aligned}
\mathbb{E}[\sup_{0\leq t\leq T}|e(t)|^{p}]\leq C\varDelta te^{CT}+2^{p}\frac{\rho p}{\varrho}C+\frac{\varrho-p}{\varrho\rho^{\frac{p}{\varrho-p}}}\frac{2C}{b^p}.	
\end{aligned}		
\end{equation}
Now, for any given $ \epsilon>0$, $\rho $ can be chosen fairly small such that 
\begin{equation}
	2^{p}\frac{\rho p}{\varrho}C<\frac{\epsilon}{3}
\end{equation}		
and $ b $ is chosen sufficiently large that
\begin{equation}
	\frac{\varrho-p}{\varrho\rho^{\frac{p}{\varrho-p}}}\frac{2C}{b^p}<\frac{\epsilon}{3}
\end{equation}	
and finally choosing $ \varDelta t $ very small such that 
\begin{equation}
C\varDelta te^{CT}<\frac{\epsilon}{3}.
\end{equation}	
Therefore
\begin{equation}
	\mathbb{E}[\sup_{0\leq t\leq T}|e(t)|^{p}]\leq \epsilon	
\end{equation}	
and the proof is complete.

\section{Almost sure exponential stability of compensated split-step theta method}
In this section, the almost sure exponential stability of the exact and numerical solutions of the stochastic pantograph model with Poisson random measure will be presented.\\\\
\textbf{Assumption 5.1.} \textit{There exist positive constants $ \xi_1 $ and $ \xi_2 $ such that for any $ p, $ $ q \in \mathbb{R}^{m}$
\begin{equation}
\label{assump4.1}
\begin{aligned}
	2\langle p,\bar{f}(p,q)\rangle+|g(p,q)|^{2}+\int_{U}|h(p,q,u)|^{2}\pi(du)\leq-2\xi_1|p|^2+2\xi_2|q|^2
\end{aligned}	
\end{equation}}
\textbf{Definition 5.1} \cite{mao2007stochastic}. The solution $ x(t) $ to the stochastic pantograph model with Poisson random measure (\ref{model}) possesses the almost sure exponential stability property if there is a positive constant $ \gamma_1 $ such that, for any initial value $ x_0 $
\begin{equation}
	\limsup_{t\rightarrow \infty}\frac{\log |x(t)|}{t}\leq -\gamma_1 \quad, a.s.
\end{equation}	
\textbf{Definition 5.2} \cite{mao2007stochastic}. The numerical solution $ \{Y_n\}_{n\geq1} $ which approximates the stochastic pantograph model with Poisson random measure (\ref{model}) possesses the almost sure exponential stability property if there is a positive constant $ \gamma_2 $ such that, for a given stepsize $ \varDelta t $ and any initial value $ x_0 $
\begin{equation}
	\limsup_{n\rightarrow \infty}\frac{\log |Y_n|}{n\varDelta t}\leq -\gamma_2 \quad a.s.
\end{equation}	
\textbf{Theorem 5.1.} \textit{ Let Assumptions 2.1 and 5.1 hold. If $ \xi_1>\xi_2\geq0 $, then for any initial data $ x_0 $, the solution $ x(t) $ to the stochastic pantograph model with Poisson random measure (\ref{model}) possesses the property of almost sure exponential stability
\begin{equation}
	\limsup_{t\rightarrow \infty}\frac{\log |x(t)|}{t}\leq -\frac{\gamma}{2} \quad a.s.
\end{equation}	
where $ \gamma=\min(1,\xi_1-\xi_2) $. }\\\\
\textbf{Proof.} The proof of the theorem can be attained by the same concept and techniques used in Theorem 3.1 of \cite{guo2018almost}.\\\\
\textbf{Theorem 5.2.} \textit{Let Assumptions 2.1 and 5.1  hold and $ \theta\in[0,\frac{1}{2}] $. Assume also, the existence of  positive constants $ \kappa_1 $ and $ \kappa_2 $ such that 
\begin{equation}
\label{asump4.2bar}
\begin{aligned}
|\bar{f}(p,q)|^2\leq\kappa_1|p|^2+\kappa_2|q|^2
\end{aligned}
\end{equation}	
for all $ p, $ $ q \in \mathbb{R}^{m}$ and $ \xi_1>\xi_2(1+[1/\eta])$. Then, there exists a constant $ \varDelta t^{*} $ such that for all $ \varDelta t \in (0,\varDelta t^*) $ and initial data $ x_0 $, the approximate solution defined by (\ref{scheme2}) possesses the property 
\begin{equation}
	\label{thm4.1}
	\limsup_{n\rightarrow \infty}\frac{\log |Y_n|}{n\varDelta t}\leq -\frac{\gamma_2}{2} \quad a.s.
\end{equation}	
where $ \gamma_2=\min(1, \zeta^{*} ) $ with $ \zeta^{*} $ satisfying
\begin{equation}
	\begin{aligned}
	\zeta^{*}(1+&\theta\varDelta t)(1+\theta\varDelta t\kappa_1)-2\xi_1+(1-2\theta)\kappa_1\varDelta t\\
	&\quad+(\zeta^{*}(1+\theta\varDelta t)\theta\varDelta t\kappa_2+2\xi_2+(1-2\theta)\kappa_2\varDelta t)(1+[1/\eta])=0
	\end{aligned}
\end{equation}
and	
\begin{equation}
	\lim_{\varDelta t\rightarrow 0}\zeta^{*}=2\xi_1-2\xi_2(1+[1/\eta]).
\end{equation}}
\textbf{Proof.} From (\ref{scheme2}), we can obtain
\begin{equation}
	\label{scheme110}
	\begin{aligned}		
		|Y_{n+1}|^{2}&=|Y_{n}|^{2}+(2\langle Y^{*}_{n}, f(Y^{*}_{n},Y^{*}_{[\eta n]})+\int_U h(Y^{*}_{n},Y^{*}_{[\eta n]},u)\pi(du)\rangle \\
		&\quad+|g(Y^{*}_{n},Y^{*}_{[\eta n]})|^{2}
		+\int_U |h(Y^{*}_{n},Y^{*}_{[\eta n]},u)|^{2}\pi(du))\varDelta t\\
		&\quad+(1-2\theta)|\bar{f}(Y^{*}_{n},Y^{*}_{[\eta n]})\varDelta t|^{2}
		+N_{n}
	\end{aligned}		
\end{equation}
where
\begin{equation}
	\begin{aligned}
		N_{n}&=|g(Y^{*}_{n},Y^{*}_{[\eta n]})|^{2}(|\varDelta W_n|^2-\varDelta t)\\
		&\quad+\int_{t_n}^{t_{n+1}}\int_{U}|h(Y^{*}_{n},Y^{*}_{[\eta n]},u)|^{2}(|\widetilde{N}(ds,du)|^{2}-\pi(du)ds)\\
		&\quad+2\langle Y_n+\bar{f}(Y^{*}_{n},Y^{*}_{[\eta n]})\varDelta t,g(Y^{*}_{n},Y^{*}_{[\eta n]})\varDelta W_n\rangle\\
		&\quad+2\langle Y_n+\bar{f}(Y^{*}_{n},Y^{*}_{[\eta n]})\varDelta t,\int_{t_n}^{t_{n+1}}\int_{U}h(Y^{*}_{n},Y^{*}_{[\eta n]},u)\widetilde{N}(ds,du)\rangle\\
		&\quad+2\langle g(Y^{*}_{n},Y^{*}_{[\eta n]})\varDelta W_n,\int_{t_n}^{t_{n+1}}\int_{U}h(Y^{*}_{n},Y^{*}_{[\eta n]},u)\widetilde{N}(ds,du)\rangle.
	\end{aligned}	
\end{equation}	
Then, by applying Assumption 5.1 and (\ref{asump4.2bar}), we get 
\begin{equation}
	\label{equality2}
	\begin{aligned}
		|Y_{n+1}|^{2}-|Y_n|^{2}&\leq-(2\xi_1-(1-2\theta)\kappa_1\varDelta t)\varDelta t|Y_n|^{2}\\
		&\quad+(2\xi_2+(1-2\theta)\kappa_2\varDelta t)\varDelta t|Y^{*}_{[\eta n]}|^{2}\\
		&\quad+N_n.
	\end{aligned}	
\end{equation}	
From (\ref{scheme1}) and (\ref{asump4.2bar}), we have
\begin{equation}
	\label{equality1}
	\begin{aligned}
		|Y_n|^{2}&\leq(1+\theta\varDelta t)|Y^{*}_{n}|^{2}+(1+\theta\varDelta t)\theta\varDelta t|\bar{f}(Y^{*}_{n},Y^{*}_{[\eta n]})|^{2}\\
		&\leq(1+\theta\varDelta t)(1+\theta\varDelta t\kappa_1)|Y^{*}_{n}|^{2}+(1+\theta\varDelta t)\theta\varDelta t \kappa_2|Y^{*}_{[\eta n]}|^{2}.
	\end{aligned}	
\end{equation}	
Then, for any constant $\zeta>0$, we have the following
\begin{equation}
	\label{inequality}
	\begin{aligned}
		e^{\zeta(n+1)\varDelta t}|Y_{n+1}|^2-e^{\zeta n\varDelta t}|Y_n|^{2}&=e^{\zeta(n+1)\varDelta t}(|Y_{n+1}|^2-|Y_n|^2)\\
		&\quad+(1-e^{-\zeta\varDelta t})e^{\zeta(n+1)\varDelta t}|Y_n|^{2}.
	\end{aligned}	
\end{equation}	
By plugging in (\ref{equality2}) and (\ref{equality1}) into (\ref{inequality}), we obtain
\begin{equation}
	\label{main40}
	\begin{aligned}
		e^{\zeta(n+1)\varDelta t}|Y_{n+1}|^2-e^{\zeta n\varDelta t}|Y_n|^{2}&\leq e^{\zeta(n+1)\varDelta t}\Biggl[A(\zeta,\varDelta t)\varDelta t|Y^{*}_{n}|^{2}\\
		&\quad+B(\zeta,\varDelta t)\varDelta t |Y^{*}_{[\eta n]}|^{2}+N_{n}\Biggl]
	\end{aligned}	
\end{equation}
where
\begin{equation}
	A(\zeta,\varDelta t)=\zeta(1+\theta\varDelta t)(1+\theta\varDelta t\kappa_1)-2\xi_1+(1-2\theta)\kappa_1\varDelta t
\end{equation}	
and
\begin{equation}
	B(\zeta,\varDelta t)=\zeta(1+\theta\varDelta t)\theta\varDelta t\kappa_2+2\xi_2+(1-2\theta)\kappa_2\varDelta t.
\end{equation}	
Taking the summation from $ j=0 $ to $ j=n-1 $ to both sides of (\ref{main40}) results in
\begin{equation}
	\label{main400}
	\begin{aligned}
		e^{\zeta n\varDelta t}|Y_{n}|^2&\leq |Y_0|^{2}+A(\zeta,\varDelta t)\varDelta t \sum_{j=0}^{n-1}e^{\zeta(j+1)\varDelta t}|Y^{*}_{j}|^{2}\\
		&\quad+B(\zeta,\varDelta t)\varDelta t \sum_{j=0}^{n-1}e^{\zeta(j+1)\varDelta t} |Y^{*}_{[\eta j]}|^{2}+N^{*}_{n}
	\end{aligned}	
\end{equation}
where $ N^{*}_{n}:= \sum_{j=0}^{n-1}e^{\zeta(j+1)\varDelta t} N_{j} $ is a martingale with $ N^{*}_{0}=0 $. Noting that for any $ l \in \{0,1,...,[\eta(n-1)]\},$  $[\eta j]=l \Leftrightarrow l\leq j< l+1 $, so we deduce that the maximum number of $ j\in \{0,1,...,n-1\} $ satisfying $ [\eta j]=l $ is $ 1+[1/\eta] $. Thus
\begin{equation}
	\label{main42}
	\begin{aligned}
		\sum_{j=0}^{n-1}&e^{\zeta(j+1)\varDelta t} |Y^{*}_{[\eta j]}|^{2}\\
		&\leq (1+[1/\eta])\Biggl(\sum_{j=0}^{n-1}e^{\zeta(j+1)\varDelta t}|Y^{*}_{j}|^{2}
		-\sum_{j=[\eta(n-1)]+1}^{n-1}e^{\zeta(j+1)\varDelta t}|Y^{*}_{j}|^{2}\Biggl).
	\end{aligned}	
\end{equation}	
By plugging in (\ref{main42}) into (\ref{main400}), we obtain
\begin{equation}
	\label{help2}
	\begin{aligned}
		e^{\zeta n\varDelta t}|Y_{n}|^2&\leq |Y_0|^{2}+T(\zeta,\varDelta t)\varDelta t \sum_{j=0}^{n-1}e^{\zeta(j+1)\varDelta t}|Y^{*}_{j}|^{2}\\
		&\quad-B(\zeta,\varDelta t)\varDelta t(1+[1/\eta]) \sum_{j=[\eta(n-1)]+1}^{n-1}e^{\zeta(j+1)\varDelta t} |Y^{*}_{j}|^{2}+N^{*}_{n}
	\end{aligned}	
\end{equation}
where $ T(\zeta,\varDelta t):=A(\zeta,\varDelta t)+B(\zeta,\varDelta t)(1+[1/\eta])$. Let
\begin{equation}
	\varDelta t^{*}=
	\left\{
	\begin{array}{ll}
		\frac{2\xi_1-2\xi_2(1-[1/\eta]}{(1-2\theta)(\kappa_1+\kappa_2(1+[1/\eta]))},   & \mbox{if } \theta \in [0,\frac{1}{2}), \\
		\infty, & \mbox{if }  \theta=\frac{1}{2}.
	\end{array}
	\right.
\end{equation}	
Then for any $ \varDelta t < \varDelta t^{*} $, 
\begin{equation}
T(0,\varDelta t)=-2(\xi_1-\xi_2(1+[1/\eta]))+(1-2\theta)\varDelta t(\kappa_1+\kappa_2 (1+[1/\eta]))<0.	
\end{equation}
It is also noted that $ \frac{d}{d\zeta}T(\zeta,\varDelta t)>0 $. So, there exists a unique constant $ \zeta^{*}>0 $ satisfying $ T(\zeta^{*},\varDelta t)=0 $. Also, we have
\begin{equation}
	\lim_{\varDelta t\rightarrow 0}T(2\xi_1-2\xi_2(1+[1/\eta],\varDelta t)=0.
\end{equation}	 
Therefore, we conclude that 
\begin{equation}
	\lim_{\varDelta t\rightarrow 0}\zeta^{*}=2\xi_1-2\xi_2(1+[1/\eta]).
\end{equation}
Using the definition of $ \gamma_2 $, if $ \zeta^{*}>1 $, then $ \gamma_2=1 $ and $ T(\gamma_2,\varDelta t)<0 $. However, if $ \zeta^{*}<1$ then $ \gamma_2=\zeta^{*} $ and $ T(\gamma_2,\varDelta t)=0.$ By applying the discrete semi-martingale convergence theorem \cite{mao2006stochastic}, for any initial data $ x_0 $, $ \limsup_{n\rightarrow \infty}[|Y_0|^{2}+N^{*}_{n}]<\infty $ which yields
\begin{equation}
	\limsup_{n\rightarrow \infty}e^{\gamma_2n\varDelta t}|Y_n|^{2}\leq\limsup_{n\rightarrow \infty}[|Y_0|^{2}+N^{*}_{n}]<\infty, \quad a.s.
\end{equation}	
Therefore, there is a finite positive random variable $ \upsilon $ such that
\begin{equation}
  \limsup_{n\rightarrow \infty}e^{\gamma_2 n\varDelta t}|Y_n|^{2}<\upsilon
\end{equation}
which implies the required assertion (\ref{thm4.1}). The proof is complete.\\\\
\textbf{Theorem 5.3.} \textit{Let Assumptions 2.1 and 5.1  hold and $ \theta\in(\frac{1}{2},1]$. Then, for any stepsize $ \varDelta t>0 $ and initial data $ x_0 $, the approximate solution defined by (\ref{scheme2}) possesses the property
	\begin{equation}
		\label{thm4.1dash}
		\limsup_{n\rightarrow \infty}\frac{\log |Y_n|}{n\varDelta t}\leq -\frac{\gamma_2}{2} \quad a.s.
	\end{equation}	
	where $ \gamma_2=\min(1, \zeta^{*} ) $ with $ \zeta^{*} $ satisfying
	\begin{equation}
			2\xi_1-2\xi_2(1+[1/\eta])-\frac{(2\theta-1)\zeta^{*}}{(2\theta-1)-\theta^{2}\zeta^{*}\varDelta t}=0
	\end{equation}
	and	
	\begin{equation}
		\lim_{\varDelta t\rightarrow 0}\zeta^{*}=2\xi_1-2\xi_2(1+[1/\eta]).
\end{equation}}
\textbf{Proof.} From (\ref{scheme2}), we can get
\begin{equation}
	\label{scheme11}
	\begin{aligned}		
		|Y_{n+1}|^{2}&=|Y_{n}|^{2}+(2\langle Y^{*}_{n}, f(Y^{*}_{n},Y^{*}_{[\eta n]})+\int_U h(Y^{*}_{n},Y^{*}_{[\eta n]},u)\pi(du)\rangle \\
		&\quad+|g(Y^{*}_{n},Y^{*}_{[\eta n]})|^{2}
		+\int_U |h(Y^{*}_{n},Y^{*}_{[\eta n]},u)|^{2}\pi(du))\varDelta t\\
		&\quad+(1-2\theta)|\bar{f}(Y^{*}_{n},Y^{*}_{[\eta n]})\varDelta t|^{2}
		+M_{n}
	\end{aligned}		
\end{equation}
where
\begin{equation}
	\begin{aligned}
		M_{n}&=|g(Y^{*}_{n},Y^{*}_{[\eta n]})|^{2}(|\varDelta W_n|^2-\varDelta t)\\
		&\quad+\int_{t_n}^{t_{n+1}}\int_{U}|h(Y^{*}_{n},Y^{*}_{[\eta n]},u)|^{2}(|\widetilde{N}(ds,du)|^{2}-\pi(du)ds)\\
		&\quad+2\langle Y_n+\bar{f}(Y^{*}_{n},Y^{*}_{[\eta n]})\varDelta t,g(Y^{*}_{n},Y^{*}_{[\eta n]})\varDelta W_n\rangle\\
		&\quad+2\langle Y_n+\bar{f}(Y^{*}_{n},Y^{*}_{[\eta n]})\varDelta t,\int_{t_n}^{t_{n+1}}\int_{U}h(Y^{*}_{n},Y^{*}_{[\eta n]},u)\widetilde{N}(ds,du)\rangle\\
		&\quad+2\langle g(Y^{*}_{n},Y^{*}_{[\eta n]})\varDelta W_n,\int_{t_n}^{t_{n+1}}\int_{U}h(Y^{*}_{n},Y^{*}_{[\eta n]},u)\widetilde{N}(ds,du)\rangle.
	\end{aligned}	
\end{equation}
We also, note that from (\ref{scheme1})
\begin{equation}
	\label{help}
	\bar{f}(Y^{*}_{n},Y^{*}_{[\eta n]})\varDelta t=\frac{1}{\theta}(Y^{*}_{n}-Y_{n}).
\end{equation}
By utilizing Assumption 5.1 and (\ref{help}), we can rewrite (\ref{scheme11}) as follows	
\begin{equation}
	\label{scheme12}
	\begin{aligned}		
		|Y_{n+1}|^{2}&\leq(1+\frac{1-2\theta}{\theta^2})|Y_{n}|^{2}-(\frac{2\theta-1}{\theta^2}+2\xi_1\varDelta t)|Y^{*}_{n}|^2 \\
		&\quad+2\xi_2\varDelta t|Y^{*}_{[\eta n]}|^{2}
		+\frac{2\theta-1}{\theta^2}2\langle Y_n,Y^{*}_{n}\rangle
		+M_{n}.
	\end{aligned}		
\end{equation}
For any constant $\zeta>0$, we have 
\begin{equation}
	\label{scheme13}
	\begin{aligned}		
		e^{\zeta(n+1)\varDelta t}&|Y_{n+1}|^{2}-e^{\zeta n\varDelta t}|Y_{n}|^{2}\\
		&\leq e^{\zeta(n+1)\varDelta t}\Biggl\{\Biggl(\frac{1-2\theta}{\theta^2}+\zeta\varDelta t\Biggl)|Y_n|^{2}
		-\Biggl(\frac{2\theta-1}{\theta^2}+2\xi_1\varDelta t\Biggl)|Y^{*}_{n}|^2 \\ 
		&\quad+2\xi_2\varDelta t|Y^{*}_{[\eta n]}|^{2}
		+\frac{2\theta-1}{\theta^2}2\langle Y_n,Y^{*}_{n}\rangle
		+M_{n}\Biggl\}.
	\end{aligned}		
\end{equation}
By summing up the above inequality (\ref{scheme13}) from $ j=0 $ to $ j=n-1 $ and proceeding by the same approach to get inequality (\ref{help2}) from inequality (\ref{main400}), we can obtain the following
\begin{equation}
	\label{scheme14}
	\begin{aligned}		
		e^{\zeta n\varDelta t}|Y_{n}|^{2}
		&\leq |Y_{0}|^{2}+(\frac{1-2\theta}{\theta^2}+\zeta\varDelta t)\sum_{j=0}^{n-1} e^{\zeta(j+1)\varDelta t}|Y_j|^{2}\\
		&\quad-(\frac{2\theta-1}{\theta^2}+(2\xi_1-2\xi_2(1+[1/\eta]))\varDelta t)\sum_{j=0}^{n-1} e^{\zeta(j+1)\varDelta t}|Y^{*}_{j}|^2 \\ 
		&\quad-2\xi_2(1+[1/\eta])\varDelta t\sum_{j=[\eta(n-1)]+1}^{n-1} e^{\zeta(j+1)\varDelta t}|Y^{*}_{j}|^2 \\
		&\quad+\frac{2\theta-1}{\theta^2}\sum_{j=0}^{n-1} e^{\zeta(j+1)\varDelta t}2\langle Y_j,Y^{*}_{j}\rangle
		+M_{n}^{*}
	\end{aligned}		
\end{equation}
where $ M_{n}^{*}=\sum_{j=0}^{n-1}e^{\zeta(j+1)\varDelta t}M_{j} $ is a martingale with $ M_{0}^{*}=0$. By choosing $ \bar{\zeta}>0 $ small enough to achieve $ \frac{1-2\theta}{\theta^2}+\bar{\zeta}\varDelta t<0 $. Therefore, for any $ \zeta \in (0,\bar{\zeta}) $,
\begin{equation}
	H(\zeta):=\frac{2\theta-1}{\theta^2}+\zeta\varDelta t>0.
\end{equation}
It is also noted that	
\begin{equation}
	\label{help3}
	\begin{aligned}
\frac{2\theta-1}{\theta^2}2\langle Y_j,Y^{*}_{j}\rangle&\leq\frac{2\theta-1}{\theta^2}\Biggl\{\frac{\theta^2	H(\zeta)}{2\theta-1}|Y_j|^{2}+\frac{2\theta-1}{\theta^2	H(\zeta)}|Y^{*}_{j}|^{2}\Biggl\}\\
&=H(\zeta) |Y_n|^{2}+\frac{(2\theta-1)^2}{\theta^4H(\zeta)}|Y^{*}_{j}|^{2}
\end{aligned}
\end{equation}	
and
\begin{equation}
	\label{help4}
	\frac{1}{\varDelta t} \Biggl(\frac{2\theta-1}{\theta^2}-\frac{(2\theta-1)^2}{\theta^4H(\zeta)}\Biggl)=-\frac{(2\theta-1)\zeta}{(2\theta-1)-\theta^2\zeta\varDelta t}.
\end{equation}	
Therefore, by utilizing (\ref{help3}) and (\ref{help4}), inequality (\ref{scheme14}) can be rewritten as 
\begin{equation}
	\label{scheme15}
	\begin{aligned}		
		e^{\zeta n\varDelta t}|Y_{n}|^{2}
		&\leq |Y_{0}|^{2}+M_{n}^{*}+\phi(\zeta,\varDelta t)\varDelta t\sum_{j=0}^{n-1} e^{\zeta(j+1)\varDelta t}|Y^{*}_{j}|^2\\
	\end{aligned}		
\end{equation}
where
\begin{equation}
	\label{help5}
	\phi(\zeta,\varDelta t)=-\Biggl(2\xi_1-2\xi_2(1+[1/\eta]) \Biggl)+\frac{(2\theta-1)\zeta}{(2\theta-1)-\theta^2\zeta\varDelta t}.
\end{equation}
From (\ref{help5}), we have
\begin{equation}
	\phi(0,\varDelta t)=-2(\xi_1-\xi_2(1+[1/\eta]))<0.
\end{equation}
It is also noted that $ \frac{d}{d\zeta}\phi(\zeta,\varDelta t)>0 $. So, there exists a unique constant $ \zeta^{*}>0 $ satisfying $ \phi(\zeta^{*},\varDelta t)=0 $. Also, we have
\begin{equation}
	\lim_{\varDelta t\rightarrow 0}\phi(2\xi_1-2\xi_2(1+[1/\eta],\varDelta t)=0
\end{equation}	 
which implies
\begin{equation}
	\lim_{\varDelta t\rightarrow 0}\zeta^{*}=2\xi_1-2\xi_2(1+[1/\eta]).
\end{equation}
Using the definition of $ \gamma_2 $, if $ \zeta^{*}>1 $, then $ \gamma_2=1 $ and $ \phi(\gamma_2,\varDelta t)<0 $. However, if $ \zeta^{*}<1$ then $ \gamma_2=\zeta^{*} $ and $ \phi(\gamma_2,\varDelta t)=0.$ By applying the discrete semi-martingale convergence theorem \cite{mao2006stochastic}, for any initial data $ x_0 $, $ \limsup_{n\rightarrow \infty}[|Y_0|^{2}+M^{*}_{n}]<\infty$, which yields
\begin{equation}
	\limsup_{n\rightarrow \infty}e^{\gamma_2n\varDelta t}|Y_n|^{2}\leq\limsup_{n\rightarrow \infty}[|Y_0|^{2}+M^{*}_{n}]<\infty, \quad a.s.
\end{equation}	
Therefore, there is a finite positive random variable $ \upsilon $ such that
\begin{equation}
	\limsup_{n\rightarrow \infty}e^{\gamma_2 n\varDelta t}|Y_n|^{2}<\upsilon
\end{equation}
which implies the required assertion (\ref{thm4.1dash}). The proof is complete.
\section{Numerical examples}
\textbf{Example 6.1.} Consider the following stochastic pantograph differential equation with Poisson random measure
\begin{equation}
	\label{example}
	\begin{aligned}
		dx(t)&=-3x(t)dt-\sin(x(\eta t))dW(t)+\int_{U}x(\eta t)\cos uN(dt,du),\quad t\in[0,2]
	\end{aligned}
\end{equation}	
with initial data $ x(0)=1$, $\eta=0.5 $ and compensator given by $ \pi(du)dt=\lambda f(u)dudt$, where $ \lambda=2 $ and $ f(u) $ is the pdf of a lognormal random variable
\begin{equation*}
	f(u)=\frac{1}{\sqrt{2\pi}u}e^{-\frac{(\ln u)^2}{2}},\quad 0\leq u<\infty.
\end{equation*}	
Because it is difficult to get the analytical solution of Equation (\ref{example}), we need to solve Equation (\ref{example}) by split-step $ \theta $ method with a sufficiently small stepsize $ (dt=1/2^{13}) $ and identify its output as the exact solution for the error comparison. Also, for stepsizes $\varDelta t=2^{q+1}dt$, where $(4\leq q\leq7)$, 2000 sample paths are used to simulate the numerical solution of Eq.(\ref{example}). For the $ i$th sample path, let the exact solution of (\ref{example}) be represented by $ Y(t_n,\omega_i) $ at time $ t=t_n $ and $ Y_n(\omega_i) $ represents the compensated split-step theta approximation at the $n$th step. Let $ \epsilon $ represents the error in the $p$th moment and by the law of large numbers, the error in the $ p$th moments at final time $T=2$ will be $ \epsilon(T)\approx\frac{1}{2000}\sum_{i=0}^{2000}|Y(T,\omega_i)-Y_N(\omega_i)|^p$. The mean-square errors of the compensated split-step theta method at the final time $ T=2 $  with $ \theta=0.5$, $ 0.75$, $ 1 $ are calculated and listed in Table 1. Also, Figure 1 depicts the plot of the mean-square errors against $ \varDelta t $ on a log-log scale and as a reference, a dashed line of slope $ 1/2 $ is added.\\
\begin{table}[!ht]
	\centering
	\caption[loftitle]{Mean-square errors for Eq.(\ref{example}).}
	\label{table1}
	\begin{tabular}{lccccc}
		\hline \\
		$ \Delta t $ & $ 2^{-5} $ & $ 2^{-6} $ & $ 2^{-7} $ & $ 2^{-8} $\\
		\hline \\
		$ \theta =0.5 $ & 0.02081807 &0.01251295 &0.00397214&0.00171497\\
		$ \theta =0.75 $ &0.03218582& 0.02134174 & 0.01089625 & 0.00316911 \\
		$ \theta =1 $ & 0.04782825 & 0.02441182 & 0.01327272 & 0.00401351 \\\\
		\hline
	\end{tabular}	
\end{table}	
\begin{figure}[!ht]
	\centering
	\begin{subfigure}[h]{0.30\columnwidth}
		\includegraphics[width=\textwidth]{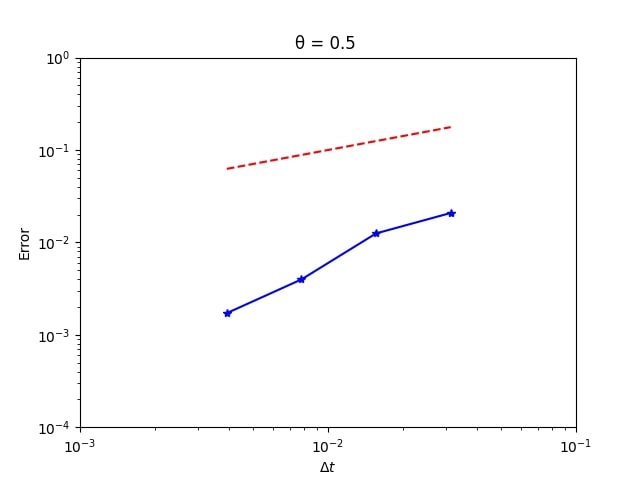}
		\caption{}
		\label{figure_theta_0.5 }
	\end{subfigure}
	~
	\begin{subfigure}[h]{0.30\columnwidth}
		\includegraphics[width=\textwidth]{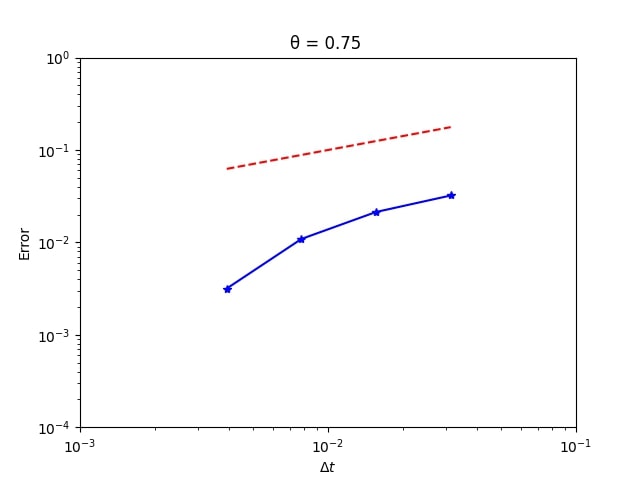}
		\caption{}
		\label{figure_theta_0.75 }
	\end{subfigure}
	~
	\begin{subfigure}[h]{0.30\columnwidth}
		\includegraphics[width=\textwidth]{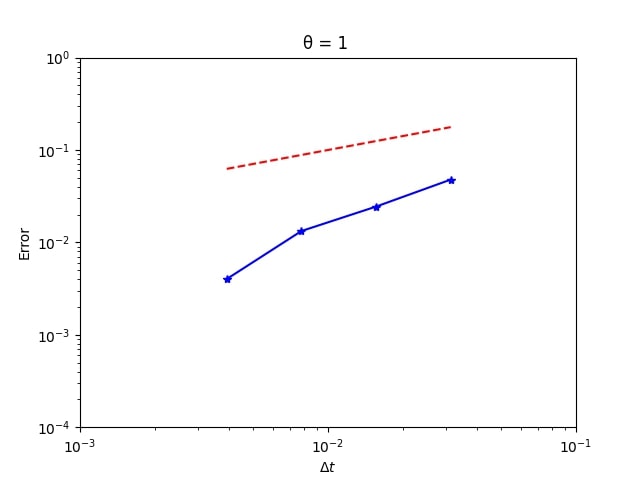}
		\caption{}
		\label{figure_theta_1 }
	\end{subfigure}
	\caption{Log-log plot of mean square errors versus stepsize $ \varDelta t $ for Eq.(\ref{example}).}
\end{figure}	
\\
\textbf{Example 6.2.}  Consider the following stochastic pantograph model with Poisson random measure
\begin{equation}
	\label{example2}
	\begin{aligned}
		dx(t)&=(-3x(t)+0.6x(\eta t))dt+x(t)\cos(x(\eta t))dW(t)\\
		&\quad+\int_{U}(0.1x(\eta t))uN(dt,du),\quad t\geq0
	\end{aligned}
\end{equation}
with initial data $ x(0)=2$, $\eta=0.5$, $ U=(0,\infty)$. The compensator is also given by $ \pi(du)dt=\lambda f(u)dudt$ where $ \lambda=2 $ and $ f(u) $ is the pdf of a lognormal random variable
\begin{equation*}
	f(u)=\frac{1}{\sqrt{2\pi}u}e^{-\frac{(\ln u)^2}{2}},\quad 0\leq u<\infty.
\end{equation*}	
According to Assumption 5.1, we have
\begin{equation*}
	\label{ha}
	\begin{aligned}
		\langle p,\bar{f}&(p,q)\rangle+\frac{1}{2}|g(p,q)|^{2}+\frac{1}{2}\int_{U}|h(p,q,u)|^{2}\pi(du)\leq-\xi_1|p|^2+\xi_2|q|^2\\
		&=p(-3p+0.6q)+\int_U(0.1\lambda pq)uf(u)du+\frac{1}{2}|p\cos q|^2\\
		&\quad+\int_{U}\frac{1}{2}\lambda|0.1q|^2u^2f(u)du\\
		&\leq -3p^2+0.01\frac{\lambda}{2} e^2q^2+(0.6+0.1\lambda \sqrt{e})pq\\
		&\quad+\frac{1}{2}|p\cos q|^2\\
		&\leq -\xi_1p^2+\xi_2q^2,
	\end{aligned}	
\end{equation*}
where $ \xi_1=2.035$ and $\xi_2=0.538$ and it is noticed that the conditions of Theorem 5.1 are fulfilled with the values of $ \xi_1 $ and $ \xi_2 $. Therefore, the solution of 
Eq.(\ref{model}) will satisfy
\begin{equation*}
	\limsup_{t\rightarrow \infty}\frac{\log |x(t)|}{t}\leq -\frac{\gamma}{2} \quad a.s.
\end{equation*}	
with $ \gamma=\min(1,\xi_1-\xi_2)=1$. Also, the conditions of Theorem 5.2 are satisfied, so the compensated split-step theta method applied to Eq.(\ref{example2}) is almost sure exponential stable when $ \theta \in [0,\frac{1}{2})$. Figures 2 and 3 depict the almost sure exponential stability of Example 6.2 by plotting three different trajectories with different selection of $ \theta $ and $ \varDelta t $.\\
\begin{figure}
	\centering
	\begin{subfigure}{0.49\textwidth}
		\includegraphics[width=\textwidth]{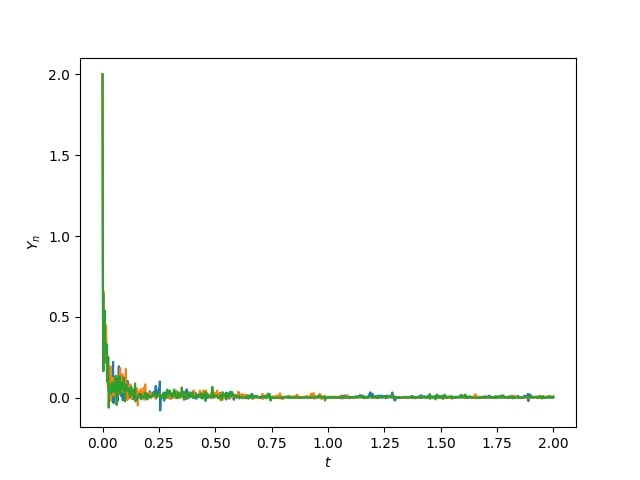}
		\caption{$\theta=0.1$, $ \varDelta t=0.2$ }
		\label{fig:first}
	\end{subfigure}
	\hfill
	\begin{subfigure}{0.49\textwidth}
		\includegraphics[width=\textwidth]{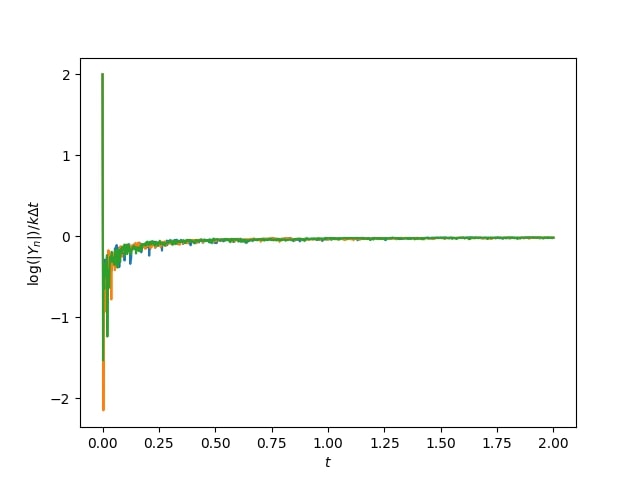}
		\caption{$\theta=0.1$, $ \varDelta t=0.2$}
		\label{fig:second}
	\end{subfigure}
	\hfill
	\begin{subfigure}{0.49\columnwidth}
		\includegraphics[width=\textwidth]{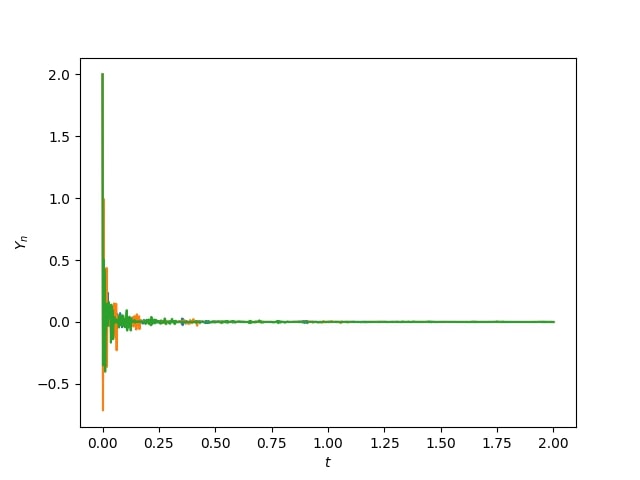}
		\caption{$\theta=0.1$, $\varDelta t=0.4$}
		\label{fig:third}
	\end{subfigure}
     \hfill
     \begin{subfigure}{0.49\columnwidth}
     	\includegraphics[width=\textwidth]{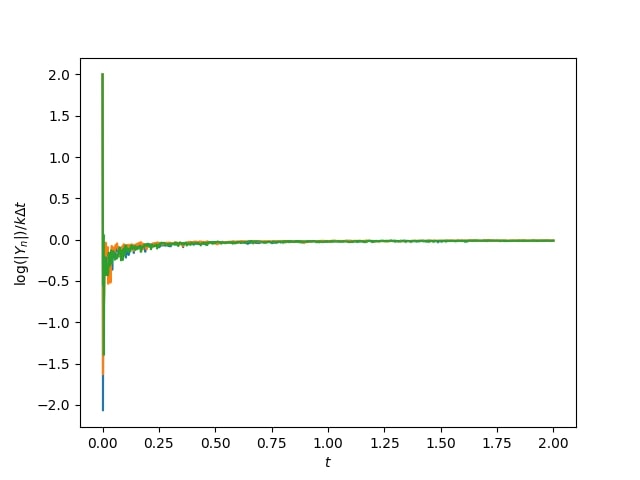}
     	\caption{$\theta=0.1$, $ \varDelta t=0.4$}
     	\label{fig:fourth}
     \end{subfigure}
	\caption{Almost sure exponential stability of (\ref{example2}) with $\theta=0.1$ and $ \varDelta t=0.2$ and $0.4$.}
	\label{fig:figures}
\end{figure}
\begin{figure}
	\centering
	\begin{subfigure}{0.49\textwidth}
		\includegraphics[width=\textwidth]{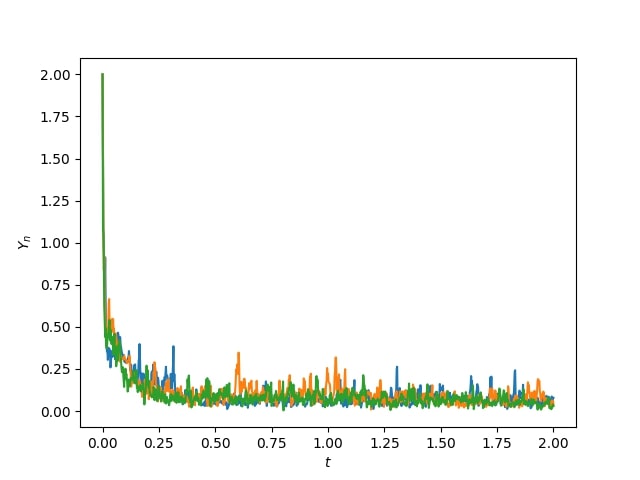}
		\caption{$\theta=0.25$, $ \varDelta t=0.2$}
		\label{fig:first1}
	\end{subfigure}
	\hfill
	\begin{subfigure}{0.49\textwidth}
		\includegraphics[width=\textwidth]{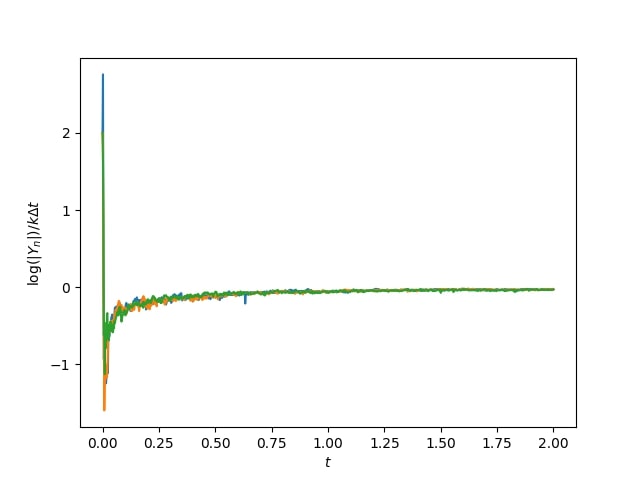}
		\caption{$\theta=0.25$, $\varDelta t=0.2$}
		\label{fig:second1}
	\end{subfigure}
	\hfill
	\begin{subfigure}{0.49\columnwidth}
		\includegraphics[width=\textwidth]{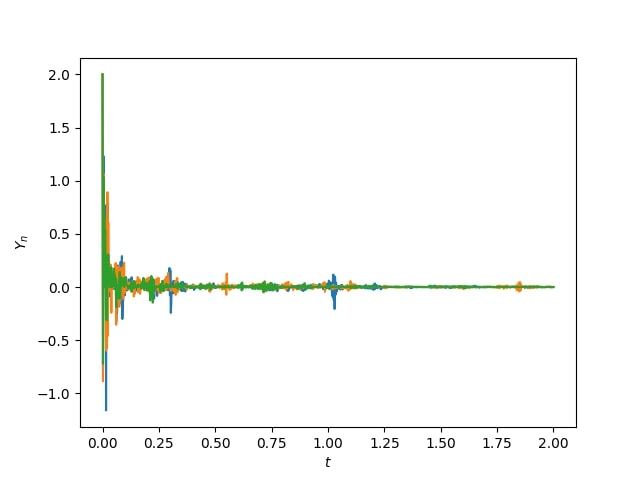}
		\caption{$\theta=0.25$, $ \varDelta t=0.4$}
		\label{fig:third1}
	\end{subfigure}
	\hfill
	\begin{subfigure}{0.49\columnwidth}
		\includegraphics[width=\textwidth]{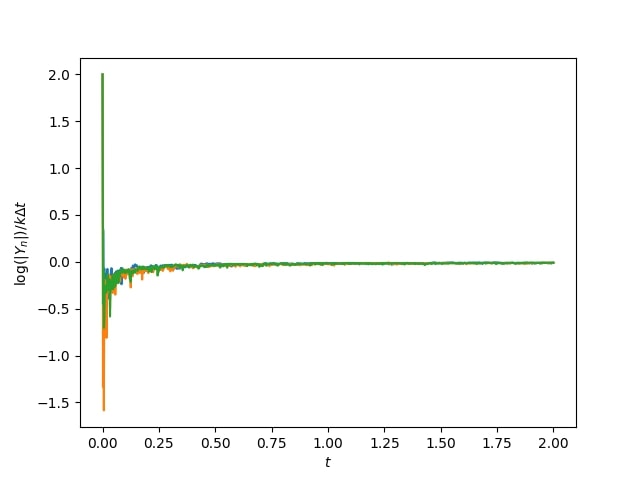}
		\caption{$\theta=0.25$, $ \varDelta t=0.4$}
		\label{fig:fourth1}
	\end{subfigure}
	\caption{Almost sure exponential stability of (\ref{example2}) with $\theta=0.25$ and $ \varDelta t=0.2$ and $0.4$.}
	\label{fig:figures1}
\end{figure}
\\
\textbf{Example 6.3.}  Consider the following stochastic pantograph model with Poisson random measure
\begin{equation}
	\label{example3}
	\begin{aligned}
		dx(t)&=(-8x(t)+2x(t)\sin^{2}(x(\eta t)))dt+\sqrt{x(t)x(\eta t)}dW(t)\\
		&\quad+\int_{U}(-0.4x(t)+0.2x(\eta t))uN(dt,du),\quad t\geq0
	\end{aligned}
\end{equation}	
with initial data $ x(0)=1$, $\eta=0.7$ and the compensator is the same as in the previous examples. By proceeding the same approach as in Example 6.2, it will be noticed that the conditions of Theorem 5.1 are fulfilled with the values of $ \xi_1=4.37 $ and $ \xi_2=1.17$. Therefore, the solution of Eq.(\ref{model}) will satisfy
\begin{equation*}
	\limsup_{t\rightarrow \infty}\frac{\log |x(t)|}{t}\leq -\frac{\gamma}{2} \quad a.s.
\end{equation*}	
with $ \gamma=\min(1,\xi_1-\xi_2)=1$. Also, the conditions of Theorem 5.3 are satisfied, so the compensated split-step theta method applied to Eq.(\ref{example3}) is almost sure exponential stable when $ \theta \in (\frac{1}{2},1]$. Figures 4 and 5 depict the almost sure exponential stability of Example 6.3 by plotting three different trajectories with different selection of $ \theta $ and $ \varDelta t $.\\
\begin{figure}[!ht]
	\centering
	\begin{subfigure}{0.49\textwidth}
		\includegraphics[width=\textwidth]{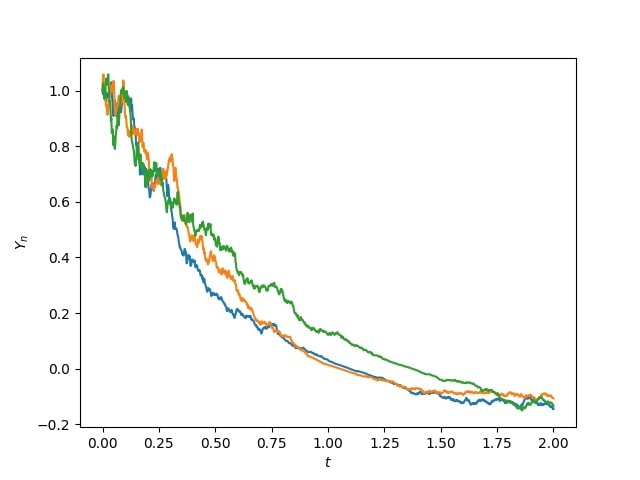}
		\caption{$\theta=0.65$, $ \varDelta t=0.7$}
		\label{fig:first11}
	\end{subfigure}
	\hfill
	\begin{subfigure}{0.49\textwidth}
		\includegraphics[width=\textwidth]{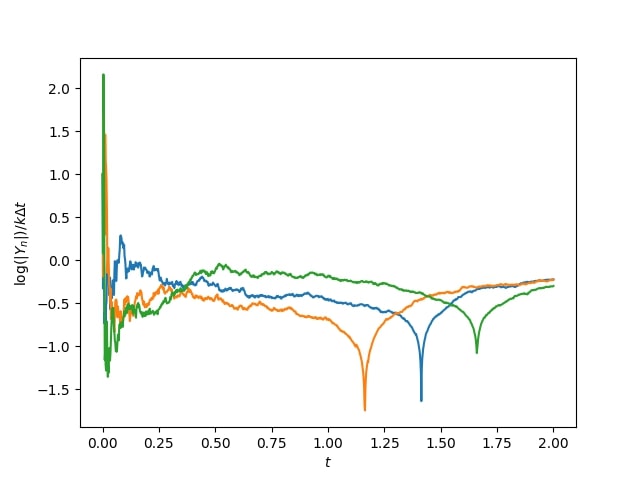}
		\caption{$\theta=0.65$, $ \varDelta t=0.7$}
		\label{fig:second11}
	\end{subfigure}
	\hfill
	\begin{subfigure}{0.49\columnwidth}
		\includegraphics[width=\textwidth]{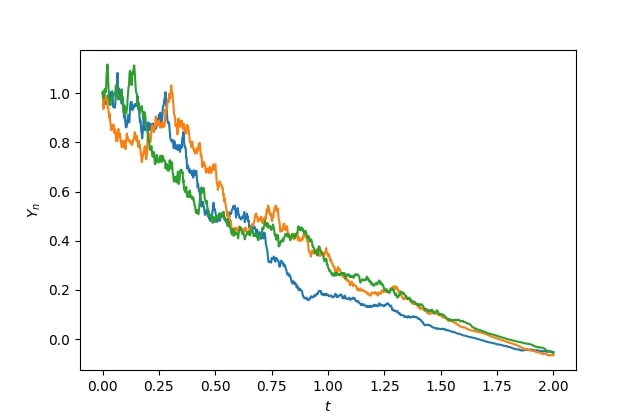}
		\caption{$\theta=0.65$, $ \varDelta t=0.9$}
		\label{fig:third11}
	\end{subfigure}
	\hfill
	\begin{subfigure}{0.49\columnwidth}
		\includegraphics[width=\textwidth]{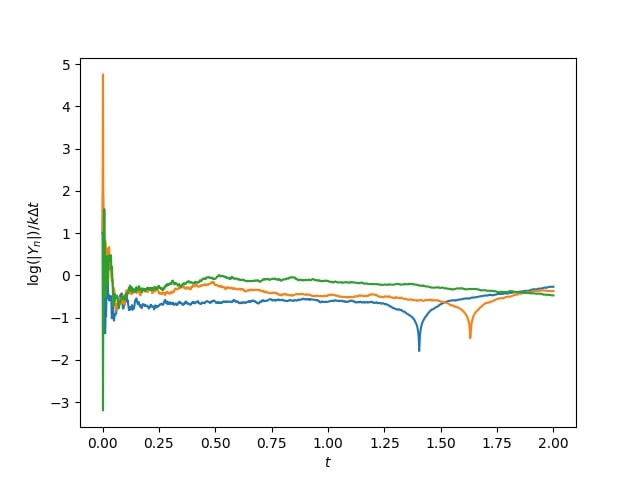}
		\caption{$\theta=0.65$, $ \varDelta t=0.9$}
		\label{fig:fourth11}
	\end{subfigure}
	\caption{Almost sure exponential stability of (\ref{example3}) with $\theta=0.65$ and $ \varDelta t=0.7$ and $0.9$.}
	\label{fig:figures11}
\end{figure}
\begin{figure}[!ht]
	\centering
	\begin{subfigure}{0.49\textwidth}
		\includegraphics[width=\textwidth]{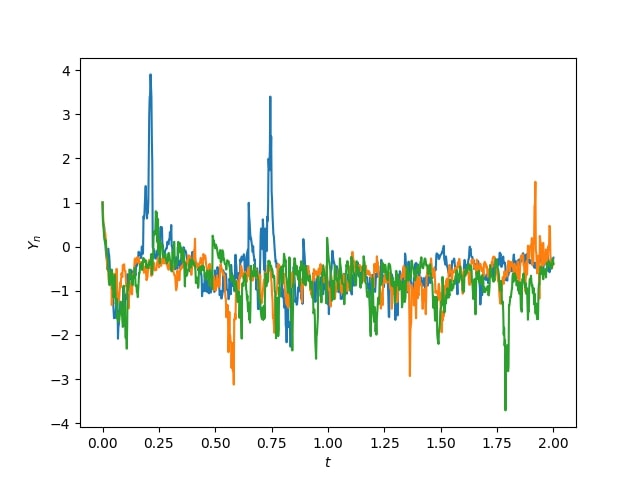}
		\caption{$\theta=0.85$, $ \varDelta t=0.7$}
		\label{fig:first111}
	\end{subfigure}
	\hfill
	\begin{subfigure}{0.49\textwidth}
		\includegraphics[width=\textwidth]{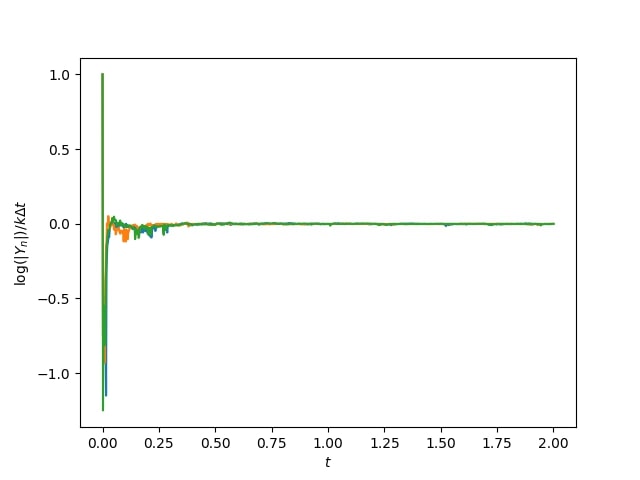}
		\caption{$\theta=0.85$, $ \varDelta t=0.7$}
		\label{fig:second111}
	\end{subfigure}
	\hfill
	\begin{subfigure}{0.49\columnwidth}
		\includegraphics[width=\textwidth]{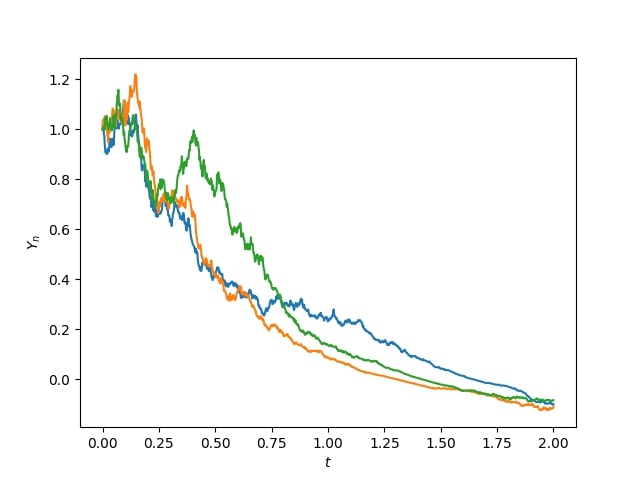}
		\caption{$\theta=0.85$, $ \varDelta t=0.9$}
		\label{fig:third111}
	\end{subfigure}
	\hfill
	\begin{subfigure}{0.49\columnwidth}
		\includegraphics[width=\textwidth]{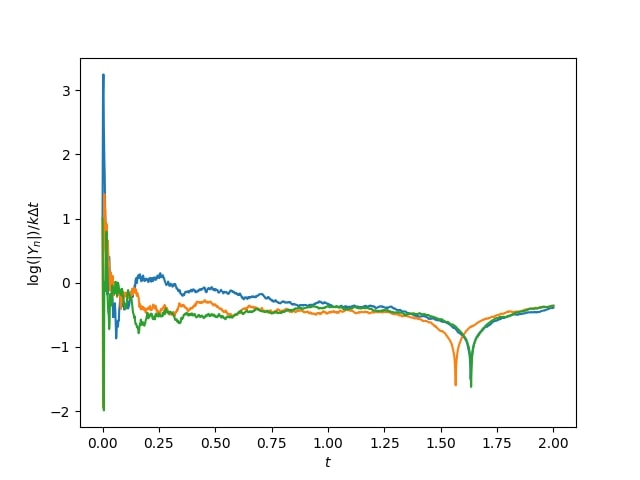}
		\caption{$\theta=0.85$, $ \varDelta t=0.9$}
		\label{fig:fourth111}
	\end{subfigure}
	\caption{Almost sure exponential stability of (\ref{example3}) with $\theta=0.85$ and $ \varDelta t=0.7$ and $0.9$.}
	\label{fig:figures111}
\end{figure}
\section{Conclusion}
In this paper, the stochastic pantograph model with Poisson random measure has been presented and the compensated split-step theta scheme was applied to it where this research contributed to the field of stochastic modelling by providing a robust and efficient numerical method for analysing stochastic pantograph models with Poisson random measure. The numerical examples demonstrated the effectiveness and practical relevance of the proposed approach, opening up new avenues for studying and understanding complex dynamical systems influenced by random factors. The numerical scheme has exhibited a non divergent attitude and manifested under the local Lipschitz condition a convergent manner in the $p$th moment to the model solution. By switching theta to either zero or one, the proposed scheme can be viewed as the compensated Euler-Maruyama method or  the compensated split-step backward Euler method respectively. Also, the sufficient conditions of almost sure exponential stability of the exact and numerical solution have been presented. By utilizing the discrete semi-martingale convergence theorem, the numerical scheme with free choice of theta can attain the almost sure exponential stability of the exact solution. However, choosing theta in range $ \theta\in[0,\frac{1}{2}] $ requires imposing extra condition on the drift side than choosing it between $ \theta\in(\frac{1}{2},1]$, for making our scheme achieving the almost sure exponential stability. Also, from the numerical section, it has been noticed that the trajectories may differ with variational choice of theta and stepsize. The conditions imposed on the numerical scheme to reproduce the almost sure exponential stability are sufficient, but not necessary ones. Therefore, our future work will focus on finding out the sufficient and necessary conditions. 
\section*{Declarations}
\textbf{Conflicts of interest}\quad The authors have not disclosed any competing interests.
\bibliographystyle{abbrv}
\bibliography{REFERENCES}
\end{document}